\documentclass{article}

\usepackage[letterpaper,nohead,top=3cm,bottom=0.5cm,left=1.5cm,right=1.5cm]{geometry}

\usepackage{amsmath,amsthm,amscd,amsfonts,amssymb,enumerate}
\usepackage{graphicx}		
\usepackage{color}
\usepackage[colorlinks]{hyperref}
\usepackage{algorithm,algorithmic}
\usepackage[titlenumbered,ruled,noend,algo2e]{algorithm2e}
\usepackage{cancel}
\usepackage{rotating}
\usepackage[graphicx]{realboxes}
\usepackage{bm}
\usepackage{tikz}
\usepackage{enumitem}
\usepackage{caption}
\usepackage{stmaryrd}
\usepackage{makecell}
\usepackage{multirow}

\usepackage{diagbox}
\usepackage{siunitx}

\usepackage{cleveref}

\usepackage[fleqn,tbtags]{mathtools}

\usepgflibrary{shapes.misc} 
\usepgflibrary[shapes.misc] 
\usetikzlibrary{shapes.misc} 
\usetikzlibrary[shapes.misc] 
\usetikzlibrary{shapes.geometric}

\newlength{\Oldarrayrulewidth}

\newlength{\defbaselineskip}
\newcommand{\setlinespacing}[1]%
           {\setlength{\baselineskip}{#1 \defbaselineskip}}

\renewcommand{\algorithmiccomment}[1]{\bgroup\hfill\{#1\}\egroup}

\newcommand\restr[2]{{
  \left.\kern-\nulldelimiterspace 
  #1 
  \littletaller 
  \right|_{#2} 
  }}
\newcommand{\littletaller}{\mathchoice{\vphantom{\big|}}{}{}{}}

\def\bkR{{\rm I\kern-.11em R}}

\theoremstyle{definition}

\begin{document}

\newcommand{\thistitle}{Improving patch selection for monolithic multigrid solvers for high-order {T}aylor-{H}ood discretizations}
\title{\thistitle}

\title{Improving patch selection for monolithic multigrid solvers for high-order {T}aylor-{H}ood discretizations}

\author {Amin Rafiei\footnote{Department of Applied Mathematics, Hakim Sabzevari University, Sabzevar, Iran. E-mail address: rafiei.am@gmail.com, a.rafiei@hsu.ac.ir, arafiei@mun.ca},
Scott MacLachlan \footnote{Department of Mathematics and Statistics, Memorial University of Newfoundland, St. John's, Canada. E-mail: s.machlachlan@mun.ca},
\\
       }
\maketitle

\begin{abstract}
Numerical simulation of incompressible fluid flows has been an active topic of research in Scientific Computing for many years, with many contributions to both discretizations and linear and nonlinear solvers.  In this work, we propose an improved relaxation scheme for higher-order Taylor-Hood discretizations of the incompressible Stokes and Navier-Stokes equations, demonstrating its efficiency within monolithic multigrid preconditioners for the linear(ized) equations.  The key to this improvement is an improved patch construction for Vanka-style relaxation introducing, for the first time, overlap in the pressure degrees of freedom within the patches.  Numerical results demonstrate significant improvement in both multigrid iterations and time-to-solution for the linear Stokes case, on both triangular and quadrilateral meshes.  For the nonlinear Navier-Stokes case, we show similar improvements, including in the number of nonlinear iterations needed in an inexact Newton method.  These improvements enable three-dimensional numerical simulations, which also presented.
\end{abstract}

%
%

\section{Introduction}
\label{sec:intro}

Numerical simulation of incompressible fluid and solid mechanics has long served as one of the main motivating applications in Scientific Computing, with the first papers published more than 50 years ago~\cite{ChorinNavStok, MarkerAndCell}.  Since then, significant advances have been made, with the development of more accurate discretization frameworks~\cite{BrezziStokes, SEMIncompFluidFlow2, JohnMixedCGDGStokes, ScottVogelius2, TaylorHoodINFSUP, CrouzeixRaviart} and efficient solvers for the resulting linear and non-linear systems of equations~\cite{DiosBrezziMariniXuZikat, ElementWiseExtendedVanka, AugLagrangeNavStok, SMAIJCMNavStok, VankaMonolithicMultigridStokes, NavierStokesVolkerTobiska, Vanka1, VankaType, VankaStokesNavierStokes, UzawaGmeinerHuberJohnRudeWohlmuth}.  In this work, we consider the important question of designing efficient preconditioners for the linear or linearized systems of equations that result from high-order (generalized) Taylor-Hood discretizations of the steady-state Stokes and Navier-Stokes equations.  

Much work has been done on the development of preconditioners for a wide variety of discretizations of the steady-state Stokes, Oseen, and Navier-Stokes equations.  Block-factorization preconditioners~\cite{Elman2014, SILVESTER2001261, doi:10.1137/040608817} are based on the idea that the Schur complements of the resulting saddle-point systems can often be approximated, either by a simple mass matrix~\cite{RVerfurth_1984a} or a corresponding convection-diffusion style operator~\cite{SILVESTER2001261, doi:10.1137/040608817}.  Combining this approach with augmented Lagrangian techniques has led to very robust preconditioners in recent years~\cite{doi:10.1137/050646421, AugLagrangeNavStok, SMAIJCMNavStok}.  A second family of effective preconditioners are monolithic multigrid methods~\cite{ABrandt_NDinar_1979a}, which rely on the definition of effective all-at-once relaxation schemes to complement coupled coarse-grid correction in a standard multigrid iteration~\cite{MultigridTutorial}.  The choice of relaxation scheme within monolithic multigrid is, of course, of utmost importance to the success of the scheme.  Commonly considered relaxation schemes include distributed relaxation~\cite{ABrandt_NDinar_1979a, DistributedRelaxBacutaVassilevskiZhang, DistributedRelaxWangChen}, Uzawa relaxation~\cite{UzawaType}, Braess-Sarazin relaxation~\cite{BraessSarazin, ZulehnerBraessSarazin}, and Vanka relaxation~\cite{ElementWiseExtendedVanka, NavierStokesVolkerTobiska, VankaMonolithicMultigridStokes, AdditiveVanka, Vanka1, VankaType, VankaStokesNavierStokes, HighOrderASMStokes}.  Here, we focus on the class of additive Vanka-style relaxation schemes~\cite{AdditiveVanka, PFarrell_etal_2019a}, extending the work from~\cite{HighOrderASMStokes} to demonstrate improved performance in the resulting preconditioners.

Vanka relaxation schemes were first proposed in the 1980's~\cite{Vanka1, VankaType} as nonlinear relaxation schemes for full approximation storage (or full approximation scheme) multigrid methods applied directly to the nonlinear equations resulting from discretizing the Navier-Stokes equations with the marker-and-cell (MAC) finite-difference discretization~\cite{MarkerAndCell}.  The fundamental idea of the relaxation scheme is to extend pointwise (weighted) Jacobi or Gauss-Seidel relaxation to a saddle-point system using the framework of overlapping Schwarz methods.  Thus, instead of inverting only the diagonal of the (linearized) discretization matrix (which may not be invertible), overlapping subsets of the degrees of freedom (DoFs) are inverted with each Vanka ``patch''.  For the MAC-scheme finite-difference discretization, these patches were defined as the cell-centred pressure, plus the four normal-direction velocity DoFs on the edges of each cell.  This results in each velocity DoF being included in two patches, while the pressure DoFs each belong to a unique patch.  This was extended to the non-conforming Crouzeix-Raviart discretization in~\cite{NavierStokesVolkerTobiska}, where a similar strategy of defining patches for each cell-centred pressure DoF with all adjacent edge-based velocity DoFs was used.  This was extended further in~\cite{JohnPatchesStokes}, where a family of ``pressure-node oriented'' relaxation schemes was proposed.  In this approach, each Vanka patch is centred on a single pressure DoF, and includes that DoF and all velocity DoFs connected to the pressure in a row of the weak divergence operator.

A topological variant of the ``pressure-node oriented'' approach was recently proposed in~\cite{HighOrderASMStokes} for higher-order discretizations, where all pressure DoFs on a common topological object (mesh element, edge, or node in a 2D mesh) and their connected velocity DoFs are relaxed simultaneously.  This leads to a \emph{composite Vanka} relaxation scheme, in which a single sweep of the composite relaxation consists of one sweep over all nodal Vanka patches plus one sweep over all edge-based Vanka patches plus one sweep over all elemental Vanka patches.  While much of the existing Vanka literature focuses on multiplicative sweeps, \cite{HighOrderASMStokes} follows the example of~\cite{PFarrell_etal_2019a, AdditiveVanka} and uses additive sweeps over all of these patches.  Numerical results there demonstrate that this leads to efficient relaxation for a variety of discretization schemes, including the Taylor-Hood discretizations considered here.  We emphasize that, to our knowledge, \cite{HighOrderASMStokes} presents the only known Vanka-style relaxation scheme that yields effective solution for high-order Taylor-Hood discretizations, although block preconditioning approaches do also exist, as in~\cite{AugLagrangeNavStok}.


In this paper, we propose a new patch construction for Vanka relaxation within monolithic multigrid preconditioners for higher-order Taylor-Hood discretizations on triangular and quadrilateral meshes.  A key difference from the work described above is that we introduce a patch construction with overlap in both the pressure and velocity DoFs.  This results in substantially fewer patches required for effective preconditioning in comparison to the composite patches proposed in~\cite{HighOrderASMStokes}, although the ``typical'' patch is somewhat larger.  Numerical results will show that the resulting relaxation scheme yields significant improvements in both iteration counts and time-to-solution for this new approach over the composite Vanka relaxation from~\cite{HighOrderASMStokes}.

The remainder of this paper is organized as follows.  \Cref{sec:mixedFEM} reviews the mixed finite-element methodology for discretizing the Stokes equations, while~\cref{sec:NavierStokes} extends this to the nonlinear time-steady Navier-Stokes equations.  \cref{sec:monolithicMG} presents the monolithic multigrid algorithm for solving the resulting linear(ized) equations, with the new patch construction discussed in~\cref{ssec:VankaStar}.  Supporting numerical results are presented in~\cref{sec:2Dnumerics,sec:3Dnumerics} for two- and three-dimensional model problems, respectively.  Conclusions and a discussion of potential future work are given in~\cref{sec:conclusion}.

\section{Mixed finite-element discretization of the Stokes equations}\label{sec:mixedFEM}
The Stokes equations describe steady, incompressible viscous flows. For a simply connected polygonal or polyhedral domain, $\Omega \subseteq \mathbb{R}^d$ for $d = 2, 3$, the Stokes equations are given by
\begin{subequations}\label{eq:Stokes}
\begin{eqnarray}
-\nabla \cdot (2\mu \bm{\varepsilon}({\bf u})) + \nabla p & = & {\bf f}\text{ in }\Omega, \label{eq:Stokes1} \\
\nabla\cdot{\bf u} & = & 0\text{ in }\Omega, \label{eq:StokesDivFree}
\end{eqnarray}
\end{subequations}
where $\mu$ is the fluid viscosity, ${\bf u}:\Omega\rightarrow \mathbb{R}^d$ is the velocity field, $p: \Omega \rightarrow \mathbb{R}$ is the pressure field and  ${\bf f}:\Omega\rightarrow \mathbb{R}^d$ is an external body force acting on the fluid. Here, $\bm{\varepsilon}({\bf u}) = \frac{1}{2}(\nabla {\bf u} + \nabla {\bf u}^T)$ is the symmetric strain-rate tensor with $\nabla {\bf u}$ defined by $(\nabla {\bf u})_{ij} = \frac{\partial u_i}{\partial x_j}$ for $1\leq i,j \leq d$.  Equation~\eqref{eq:Stokes1} is known as the momentum balance equation, while~\eqref{eq:StokesDivFree} is known as the continuity equation, or as the incompressibility equation. No-slip Dirichlet boundary conditions are formulated as
\begin{eqnarray*}\label{eq:StokesBCforCGFEM}
{\bf u} & = & {\bf 0}\text{ on }\partial \Omega,  \label{eq:StokesDirichletNoSlip}
\end{eqnarray*}
where $\partial \Omega$ denotes the boundary of the domain. 

A mixed finite-element discretization of~\eqref{eq:Stokes} arises by choosing suitable spaces for the velocity and pressure functions and passing from the strong form above into the weak form of the equations.
Let ${\bf H^1}(\Omega) = [H^1(\Omega)]^d$, and consider the two Hilbert spaces 
\begin{eqnarray*}\label{eq:H1L02modR}
 {\bf H_0^1}(\Omega) & := & \{ {\bf v}\in {\bf H^1}(\Omega)~|~{\bf v} = {\bf 0}\text{ on }\partial \Omega\}, \\ \nonumber
L_0^2(\Omega) = L^2(\Omega) / \mathbb{R} & := & \left\{q\in L^2(\Omega)~\middle|~ \int_{\Omega}q\,dx = 0\right\},
\end{eqnarray*}
with their associated norms $\|\cdot\|_{{\bf H^1}(\Omega)}$ and $\|\cdot\|_{L^2(\Omega)}$.  Multiplying the two equations in~\eqref{eq:Stokes} by suitable test functions and integrating by parts leads to the bilinear forms
\begin{subequations}\label{eq:sdsH01}
\begin{align}
a({\bf u},{\bf v}) & = 2\mu \int_{\Omega}\bm{\varepsilon}({\bf v}):\bm{\varepsilon}({\bf u})~dx, \label{eq:bilinear_form_a}\\ 
b({\bf v},p) & = - \int_{\Omega} p(\nabla \cdot {\bf v})~dx,
\end{align}
\end{subequations}
where $a({\bf u},{\bf v})$ depends on the tensor contraction ${\bf U}:{\bf V} = \sum_{i,j=1}^d U_{ij}V_{ij}$, where ${\bf U}$ and ${\bf V}$ are $d\times d$ tensors.  With these, we can express the weak form of~\eqref{eq:Stokes} as finding ${\bf u}\in {\bf H_0^1}(\Omega)$, $p\in L_0^2(\Omega)$ such that
\begin{subequations}\label{eq:stokes_weak}
\begin{align}
a({\bf u},{\bf v}) + b({\bf v},p) & = \langle {\bf f},{\bf v} \rangle, \\
b({\bf u},q) & = 0, 
\end{align}
\end{subequations}
for all ${\bf v}\in {\bf H_0^1}(\Omega)$ and $q\in L_0^2(\Omega)$.  Well posedness of~\eqref{eq:stokes_weak} is guaranteed by standard theory (see, for example, \cite{MixedFEMAppl}), relying on the continuity of $a$ and $b$, the coercivity of $a$, and that $b$ satisfies an inf-sup condition.

To discretize~\eqref{eq:stokes_weak}, we introduce a mesh, $\mathcal{T}_h$, of $\Omega$.  In $\mathbb{R}^2$, we will consider triangular and quadrilateral meshes of polygonal domains.  In $\mathbb{R}^3$, we will focus on tetrahedral meshes of polyhedral domains, but the same approach naturally extends to hexahedral elements.  Taking ${\bf V}_h$ and $W_h$ to be ${\bf H_0^1}$- and $L_0^2$-conforming finite-element spaces over $\mathcal{T}_h$, respectively, the discrete variational form is to find $({\bf u}_h,p_h)\in {\bf V}_h \times W_h$ such that
\begin{subequations}\label{eq:VariationalDiscStokes}
\begin{eqnarray}
a({\bf u}_h,{\bf v}_h) + b({\bf v}_h,p_h) & = & \langle{\bf f},{\bf v}_h \rangle, \label{eq:DiscMomentum} \\
           b({\bf u}_h,q_h) & = & 0, \label{eq:DiscCompressibility}
\end{eqnarray}
\end{subequations}
for all ${\bf v}_h\in {\bf V}_h$ and $q_h\in W_h$.  For well-posedness of a conforming discretization, the continuity of $a$ and $b$ and coercivity of $a$ are inherited from the continuum, but we must separately prove the discrete inf-sup condition,
\begin{eqnarray}\label{eq:DiscINFSUP}
\adjustlimits\inf_{q_h\in W_h\backslash\{0\}} \sup_{{\bf v}_h\in {\bf V}_h \backslash\{\bf 0\}} \frac{b({\bf v}_h,q_h)}{\|{\bf v}_h\|_{{\bf H^1}(\Omega)}\|q_h\|_{L^2(\Omega)}}\ge \beta_h>0. 
\end{eqnarray}


When $\mathcal{T}_h$ is a simplicial mesh, we let $P_k = P_k(\mathcal{T}_h)$ be the continuous Lagrange finite-element space of order $k$ over $\mathcal{T}_h$.  Similarly, when $\mathcal{T}_h$ is a quadrilateral (or hexahedral) mesh, we let $Q_k = Q_k(\mathcal{T}_h)$ be the corresponding continuous Lagrange finite-element space of order $k$ over $\mathcal{T}_h$.  We denote their vector-valued counterparts by ${\bf P_k} = [P_k]^d$ and ${\bf Q_k} = [Q_k]^d$.  The (generalized) Taylor-Hood finite-element pairs are given by $({\bf V}_h,W_h) = ({\bf Q_k},Q_{k-1})$ on quadrilateral/hexahedral meshes and $({\bf V}_h,W_h) = ({\bf P_k},P_{k-1})$ on triangular/tetrahedral meshes, for $k\ge 2$.  These are the most common ${\bf H_0^1}$- and $L_0^2$-conforming finite-element schemes for discretizing~\eqref{eq:stokes_weak} and are known to satisfy the discrete inf-sup condition~\eqref{eq:DiscINFSUP} with $\beta_h\ge \beta_0$ as $h\rightarrow 0^+$ for any order, $k$~\cite{MixedFEMAppl}.  

Defining the velocity and pressure finite-element basis functions as 
\begin{eqnarray}\label{eq:FEMbasisVelocityPressure}
{\bf V}_h = \text{span}\{ {\bm \Phi}_1,\cdots {\bm \Phi}_n \}, & \qquad & W_h = \text{span}\{ \Psi_1,\cdots \Psi_m \},
\end{eqnarray}
the finite-element representations of the approximate solutions, ${\bf u}_h$ and $p_h$, can be written as
\begin{eqnarray*}\label{eq:FEMrepresentVelocityPressure}
{\bf u}_h = \sum_{j=1}^{n} {u}_j{\bm \Phi}_j, & \qquad & p_h = \sum_{j=1}^{m} {p}_j\Psi_j,
\end{eqnarray*}
where $\left\{u_j\right\}$ and $\left\{p_j\right\}$ are the coefficients of the solution in these bases.  We note that it is a common practice for velocities in ${\bf H_0^1}(\Omega)$ to omit basis functions associated with boundary nodes of the mesh (since the solution is zero there anyway).  However, it is common to include all pressure DoFs in the basis for $W_h$, and recover the solution $p_h \in L_0^2(\Omega)$ by post-processing.  With these, the Taylor-Hood finite-element discretization of~\eqref{eq:stokes_weak} leads to solving the symmetric indefinite saddle-point system of equations written as 
\begin{eqnarray}\label{eq:StokesSaddlePointSys}
\underbrace{\left( \begin{array}{cc}
A & B^T \\
B & 0
\end{array}\right)}_{\mathcal{A}}
\underbrace{\left( \begin{array}{c}
{\bf {u}} \\
{\bf {p}}
\end{array}\right)}_{{\bf {x}}} = {\bf {b}},
\end{eqnarray}
where $\mathcal{A}\in \mathbb{R}^{(n+m) \times (n+m)}$, ${\bf {u}} = [{u}_1,\cdots,{u}_n]^T$ and ${\bf {p}} = [{p}_1,\cdots,{p}_m]^T$. The symmetric positive-definite vector-Laplacian matrix $A = [a_{ij}]\in \mathbb{R}^{n\times n}$ and rectangular divergence matrix $B = [b_{ij}]\in \mathbb{R}^{m\times n}$ are given by 
\begin{eqnarray}\label{eq:VectorLaplacianDivergenceMats}
a_{ij} = 2\mu \int_{\Omega} \bm{\varepsilon}({\bm \Phi}_i):\bm{\varepsilon}({\bm \Phi}_j)~dx, & \qquad &  b_{ij} = - \int_{\Omega} \Psi_i \nabla \cdot {\bm \Phi}_j~dx.
\end{eqnarray}
We note that, if the normalization for $L_0^2(\Omega)$ has not been dealt with by choice of basis, then $B^T$ has a one-dimensional nullspace spanned by the constant vector.  In~\cref{sec:monolithicMG}, we develop a monolithic multigrid preconditioner for systems like those in~\eqref{eq:StokesSaddlePointSys}. 

\section{Mixed finite-element discretization of the Navier-Stokes equations}\label{sec:NavierStokes}

The nondimensionalized (cf.~\cite{ChorinNavStok}) stationary incompressible Navier-Stokes equations are given by
\begin{subequations}\label{eq:NavierStokes}
\begin{eqnarray}
-\nabla \cdot (2 \bm{\varepsilon}({\bf u})) + Re({\bf u}\cdot\nabla){\bf u} + \nabla p & = & {\bf f},\text{ in }\Omega,\\
\nabla\cdot{\bf u} & = & 0,\text{ in }\Omega,\\
{\bf u} & = & {\bf g},\text{ on }\partial \Omega,
\end{eqnarray}
\end{subequations}
where $\Omega$, ${\bf u}$, $p$, and ${\bf f}$ are defined as in the Stokes case. Here, we introduce $Re$ as the Reynolds number and ${\bf g}$ as a (potentially non-homogeneous) Dirichlet boundary condition satisfying the compatibility condition 
\begin{eqnarray*}
\int_{\partial \Omega} {\bf g}\cdot{\bf n}~ds = 0,
\end{eqnarray*}
where ${\bf n}$ is the outward-pointing unit normal vector on $\partial\Omega$.  While much of the original numerical work on the Navier-Stokes equations (including~\cite{ChorinNavStok}) makes use of the vector Laplacian in place of $\nabla \cdot (2 \bm{\varepsilon}({\bf u}))$, we keep the stress-tensor form for consistency with many recent approaches.
The Reynolds number is defined as the ratio of inertial forces to viscous forces within a fluid, defined as:
\begin{eqnarray*}
Re = \frac{\rho U L}{\mu},
\end{eqnarray*}
where
\begin{itemize}
\item $\rho$ is the density of the fluid (with units $kg/m^3$)
\item $U$ is a characteristic flow speed (with units $m/s$)
\item $L$ is a characteristic linear dimension or characteristic length (with units $m$)
\item $\mu$ is the dynamic viscosity of the fluid (with units $kg/(m\cdot s)$). 
\end{itemize}
The Reynolds number quantifies the relative importance of inertial and viscous forces for given flow conditions, and laminar flow occurs at low Reynolds numbers, where viscous forces are dominant, while turbulent flow occurs at high Reynolds numbers, where inertial forces are dominant.

For ${\bf H_0^1}$- and $L_0^2$-conforming finite-element spaces, ${\bf V}_h$ and $W_h$, the nonlinear discretized weak variational form of~\eqref{eq:NavierStokes} is to find approximate solution $({\bf u}_h,p_h)\in ({\bf V}_h,W_h)$ such that
\begin{subequations}\label{eq:VariationalDiscNavStokes}
\begin{eqnarray}
{a}({\bf u}_h,{\bf v}_h) + b({\bf v}_h,p_h) + c({\bf u}_h,{\bf u}_h,{\bf v}_h) & = & \langle{\bf f},{\bf v}_h\rangle, \label{eq:DiscMomentumNS}\\
           b({\bf u}_h,q_h) & = & 0, \label{eq:DiscCompressibilityNS}
\end{eqnarray}
\end{subequations}
for all ${\bf v}_h\in {\bf V}_h$ and $q_h\in W_h$, where we reuse the bilinear form $a(\cdot,\cdot)$ from~\eqref{eq:bilinear_form_a} but with $\mu = 1$ (since we use $Re$ to account for the viscosity in this setting), and define the trilinear form for ${\bf u}_h$, ${\bf w}_h$, and ${\bf v}_h \in {\bf V}_h$:
\begin{equation*}\label{eq:3LinearH01}
c({\bf u}_h,{\bf w}_h,{\bf v}_h) =  Re  \int_{\Omega} \left[{\bf u}_h\cdot \nabla {\bf w}_h\right]\cdot {\bf v}_h~dx.
\end{equation*}
A common iterative approach for solving \eqref{eq:VariationalDiscNavStokes} is the Newton linearization technique, which builds a sequence of approximate solutions $({\bf u}^{(k)}_h,p^{(k)}_h)\in ({\bf V}_h,W_h)$ from a given initial guess $({\bf u}^{(0)}_h,p^{(0)}_h)$. At step $k$ of this process, the nonlinear residual is given by
\begin{eqnarray*}\label{eq:NonLinResNavStokes}
R_{{\bf u}}({\bf v}_h) & = & \langle{\bf f},{\bf v}_h\rangle - {a}({\bf u}^{(k)}_h,{\bf v}_h) - b({\bf v}_h,p_h^{(k)}) - c({\bf u}^{(k)}_h,{\bf u}^{(k)}_h,{\bf v}_h), \nonumber\\
R_{p}(q_h)  & = & - b({\bf u}^{(k)}_h,q_h), 
\end{eqnarray*}
for all ${\bf v}_h\in {\bf V}_h$ and $q_h\in W_h$. The next approximate pair is
\begin{eqnarray*}\label{eq:kplus1stNewtonNS}
({\bf u}^{(k+1)}_h,p^{(k+1)}_h) = ({\bf u}^{(k)}_h + \omega_{{\bf u}}\delta{\bf u}_h,p^{(k)}_h + \omega_{p}\delta p_h), 
\end{eqnarray*}
where the two weight parameters, $\omega_{{\bf u}}$ and $\omega_{p}$, are included when Newton's method is combined with a line search~\cite{LineSearch} or other damping.  
Linearizing~\eqref{eq:VariationalDiscNavStokes} leads to solving the linear problem
\begin{subequations}\label{eq:LinearizedNewtonNS}
{\small
\begin{eqnarray}
{a}(\delta{\bf u}_h,{\bf v}_h) + b({\bf v}_h,\delta p_h) + c(\delta{\bf u}_h,{\bf u}^{(k)}_h,{\bf v}_h) + c({\bf u}^{(k)}_h,\delta{\bf u}_h,{\bf v}_h) & = R_{{\bf u}}({\bf v}_h),   \\
b(\delta{\bf u}_h,q_h) & = R_{p}(q_h), 
\end{eqnarray}
}
\end{subequations}
for the increments, $\delta{\bf u}_h$ and $\delta p_h$, for all ${\bf v}_h\in {\bf V}_h$ and $q_h\in W_h$. Following~\eqref{eq:FEMbasisVelocityPressure}, the increments can be written as
\begin{eqnarray*}\label{eq:FEMrepresentDeltaVelocityDeltaPressure}
\delta{\bf u}_h = \sum_{j=1}^{n} \Delta{u}_j{\bm \Phi}_j, & \qquad & \delta p_h = \sum_{k=1}^{m} \Delta{p}_k\Psi_k.
\end{eqnarray*}
The Taylor-Hood finite-element discretization of~\eqref{eq:LinearizedNewtonNS} then leads to the system
\begin{eqnarray}\label{eq:DiscNavStokesSaddlePointSys}
\underbrace{\left( \begin{array}{cc}
C & B^T \\
B & 0
\end{array}\right)}_{\tilde{\mathcal{A}}}
\underbrace{\left( \begin{array}{c}
\bm{\Delta{u}} \\
\bm{\Delta{p}}
\end{array}\right)}_{\tilde{\bf x}} = \tilde{\bf b},
\end{eqnarray}
where nonsymmetric matrix $\tilde{\mathcal{A}}\in \mathbb{R}^{(n+m) \times (n+m)}$, $C = {A} + N + W$, $\bm{\Delta{u}} = [\Delta{u}_1,\cdots,\Delta{u}_n]^T$, ${\bf \Delta{p}} = [\Delta{p}_1,\cdots,\Delta{p}_m]^T$ and matrices $A$ and $B$ have the same definition as in \eqref{eq:VectorLaplacianDivergenceMats}.
The vector-convection matrix $N = [n_{ij}]\in \mathbb{R}^{n \times n}$ and the Newton derivative matrix $W = [w_{ij}]\in \mathbb{R}^{n \times n}$ are given by
\begin{eqnarray*}\label{eq:VectorLaplacianConvectionDerivMats}
n_{ij} = Re \int_{\Omega} ({\bf u}^{(k)}_h \cdot \nabla{\bm \Phi}_j) \cdot {\bm \Phi}_i~dx, & \qquad & w_{ij} = Re \int_{\Omega} ({\bm \Phi}_j \cdot \nabla{\bf u}^{(k)}_h)\cdot {\bm \Phi}_i~dx.
\end{eqnarray*}
In an exact Newton method, system~\eqref{eq:DiscNavStokesSaddlePointSys} is solved by a direct method. In contrast, inexact Newton methods allow iterative approximation of the solution of this system with the stopping criterion controlled by a dynamic forcing term~\cite{InexactNewton}. The approach of Eisenstat and Walker provides flexible options for effectively choosing the forcing terms~\cite{EisenStatWalker}, which we use here. 

In this paper, we consider Newton-Krylov-Multigrid (NKM) solvers for \eqref{eq:VariationalDiscNavStokes}, using Newton's method to linearize the nonlinear system, with the Eisenstat-Walker stopping criterion used for the inner Krylov iteration, which is preconditioned by multigrid.
The strong coupling between the velocity and pressure variables  motivates us to design a monolithic multigrid preconditioner for the Stokes and Navier-Stokes equations. In the next section, we review the basic concepts of a patch-based monolithic multigrid algorithm as a solver.

\section{Monolithic patch-based geometric multigrid}\label{sec:monolithicMG}
We consider the linear systems of equations that arise from approximating solutions of either the discretized Stokes equations~\eqref{eq:VariationalDiscStokes} or Newton linearization of the discretized Navier-Stokes equations \eqref{eq:LinearizedNewtonNS} with finite-element spaces $({\bf V}_h,W_h) = ({\bf P_k},P_{k-1})$ or $({\bf V}_h,W_h) = ({\bf Q_k},Q_{k-1})$, for $k\ge 2$.  While block-factorization preconditioners have also been considered for these problems~\cite{AugLagrangeNavStok, SILVESTER2001261, doi:10.1137/040608817}, we focus here on the monolithic multigrid methodology, where we directly apply multigrid to the coupled velocity-pressure system.  In this section, we present the components of the geometric multigrid preconditioner that we develop.  In all cases, we consider only a standard multigrid V-cycle to define the preconditioner.

For a fixed polynomial order $k$, consider $({\bf V}_{h,1},W_{h,1})$ as finite-element spaces associated with the coarsest level of multigrid $V$-cycle, where the corresponding mesh, $\mathcal{T}_1$, is the finite-element mesh for the domain $\Omega$ associated with this pair. Following standard geometric $h$-refinement, we build a family of nested finite-element meshes, $\mathcal{T}_1\subset \mathcal{T}_2\subset \cdots \subset \mathcal{T}_L$, for domain $\Omega$, where $\mathcal{T}_L$ is the finest mesh of the multigrid $V$-cycle. 
The nested finite-element meshes induce nested pairs of finite-element spaces 
\begin{eqnarray}\label{eq:NestedFEMSpaces}
({\bf V}_{h,\ell -1},W_{h,\ell -1}) \subset ({\bf V}_{h,\ell },W_{h,\ell}) & \quad & \ell = 2, \cdots, L,
\end{eqnarray}
where $({\bf V}_{h,\ell },W_{h,\ell})$ is associated with mesh $\mathcal{T}_{\ell}$, at level $\ell$ of the multigrid hierarchy. 
Discretization of the weak variational forms~\eqref{eq:VariationalDiscStokes} or~\eqref{eq:LinearizedNewtonNS} at level $\ell$ of the multigrid hierarchy using the Taylor-Hood finite-element spaces $({\bf V}_{h,\ell},W_{h,\ell})$, leads to the linear systems of saddle-point equations 
\begin{eqnarray*}\label{eq:StokesSaddlePointSysMG}
\underbrace{\left( \begin{array}{cc}
A_{\ell} & B_{\ell}^T \\
B_{\ell} & 0
\end{array}\right)}_{\mathcal{A}_{\ell}}
\underbrace{\left( \begin{array}{c}
{\bf {u}}_{\ell} \\
{\bf {p}}_{\ell}
\end{array}\right)}_{x_{\ell}} = {\bf b}_{\ell}, & \quad & 
\underbrace{\left( \begin{array}{cc}
C_{\ell} & B_{\ell}^T \\
B_{\ell} & 0
\end{array}\right)}_{\tilde{\mathcal{A}}_{\ell} }
\underbrace{\left( \begin{array}{c}
\bm{\Delta{u}}_{\ell} \\
\bm{\Delta{p}}_{\ell}
\end{array}\right)}_{\tilde{\bf x}_{\ell}} = \tilde{\bf b}_{\ell},
\end{eqnarray*}
where $\mathcal{A}_{\ell}, \tilde{\mathcal{A}}_{\ell}\in \mathbb{R}^{(n_{\ell}+m_{\ell}) \times (n_{\ell}+m_{\ell})}$ and $n_{\ell}$ and $m_{\ell}$ are the dimensions of finite-element spaces ${\bf V}_{h,\ell}$ and $W_{h,\ell}$, respectively.

There are three main components of the monolithic multigrid method on each level. For the prolongation of vectors from  level $\ell-1$ to level $\ell$ of the multigrid hierarchy, we use the canonical coupled prolongation operator
\begin{eqnarray*}\label{eq:ProlongMonolithicStokes}
\mathcal{P}_{\ell} = \left(\begin{array}{cc}
P_{\bf V} & 0\\
0            & P_W
\end{array}\right),
\end{eqnarray*}
where $P_{\bf V}\in \mathbb{R}^{n_{\ell} \times n_{{\ell}-1}}$ is the matrix representation of the finite-element interpolation operator associated with the natural embedding ${\bf V}_{h,{\ell}-1}\subset {\bf V}_{h,\ell}$ and $P_W\in \mathbb{R}^{m_{\ell} \times m_{{\ell}-1}}$ is the matrix representation of the finite-element interpolation operator corresponding to the natural embedding $W_{h,{\ell}-1}\subset W_{h,\ell}$.  
The restriction operator $\mathcal{R}_{\ell}$ from level ${\ell}$ to level ${\ell}-1$  is defined as
$\mathcal{R}_{\ell} = \mathcal{P}_{\ell}^T$.  In this work, we use rediscretization to define the matrices in~\eqref{eq:StokesSaddlePointSysMG} for the operators on each level (noting that this is equivalent to a Galerkin coarsening if suitable quadrature is used for the trilinear form).  Finally, at each level of the multigrid $V$-cycle, we define relaxation by applying a fixed number of steps of FGMRES using a Vanka-style relaxation scheme (outlined in the following subsections) as a preconditioner, noting that we use FGMRES to avoid parameter choice that could otherwise be accomplished using local Fourier analysis~\cite{PFarrell_etal_2019a, JBrown_etal_2019a}.
\subsection{Monolithic patch-based relaxation schemes}

A space decomposition of the product finite-element space ${\bf V}_h\times W_h$ is given by writing
\begin{eqnarray*}\label{eq:spacedecomp}
{\bf V}_h\times W_h = \sum\limits_{i=1}^{I} \left({\bf V}_h^{(i)}\times W_h^{(i)}\right),
\end{eqnarray*}
meaning that every $({\bf v}_h,q_h)\in {\bf V}_h\times W_h$ has a (not necessarily unique) representation $({\bf v}_h,q_h) = \sum\limits_{i = 1}^{I}({\bf v}_h^{(i)},q_h^{(i)})$, for ${\bf v}_h^{(i)}\in {\bf V}_h^{(i)}$ and $q_h^{(i)}\in W_h^{(i)}$ \cite{Spacedecomp}.  Both additive and multiplicative iterative methods can be defined once the space decomposition is specified, and we focus here on the resulting additive subspace correction (or additive Schwarz) algorithm.

There are two natural ways to view these subspace correction algorithms.  The first is from the finite-element perspective.  Here, if $({\bf u}_h,p_h) \in {\bf V}_h \times W_h$ is the current approximation to the solution of the discretized Stokes equations in~\eqref{eq:VariationalDiscStokes}, then, for each subdomain, $i$, we solve for the update $({\bf e}_{\bf u}^{(i)}, e_p^{(i)}) \in {\bf V}_h^{(i)} \times W_h^{(i)}$ that satisfies
\begin{subequations}\label{eq:UpdateStokes}
\begin{eqnarray}
{a}({\bf u}_h + {\bf e}_{\bf u}^{(i)},{\bf v}_h) + b({\bf v}_h,p_h + {e}_{p}^{(i)}) & = & \langle{\bf f},{\bf v}_h\rangle, \\
           b({\bf u}_h + {\bf e}_{\bf u}^{(i)},q_h) & = & 0, 
\end{eqnarray}
\end{subequations}
for all $({\bf v}_h,q_h) \in {\bf V}_h^{(i)} \times W_h^{(i)}$.  We then update the current approximation as
\[
({\bf u}_h, p_h) \leftarrow \left({\bf u}_h + \sum_{i=1}^I w_i {\bf e}_{\bf u}^{(i)}, p_h + \sum_{i=1}^I \gamma_i e_p^{(i)}\right),
\]
where $w_i: {\bf V}_h^{(i)} \rightarrow {\bf V}_h^{(i)}$ and $\gamma_i: W_h^{(i)} \rightarrow W_h^{(i)}$ are weighting operators that allow us to apply, for example, partition of unity scaling to the resulting corrections when only some DoFs are contained in multiple subspaces within the decomposition.  We note that we do not use such scalings in the numerical results to follow, but include them here for completeness.  For the Navier-Stokes case, a similar restriction is made in~\eqref{eq:LinearizedNewtonNS}, but now replacing the Newton updates $\delta{\bf u}_h$ and $\delta p_h$ by their locally updated forms, $\delta{\bf u}_h + {\bf e}_{\bf u}^{(i)}$ and $\delta p_h + e_p^{(i)}$, for $({\bf e}_{\bf u}^{(i)}, e_p^{(i)}) \in {\bf V}_h^{(i)} \times W_h^{(i)}$ and then restricting the linearized variational form to considering $({\bf v}_h, q_h) \in {\bf V}_h^{(i)} \times W_h^{(i)}$.  The update to the Newton updates is then computed accordingly, as
\[
(\delta{\bf u}_h, \delta p_h) \leftarrow \left(\delta {\bf u}_h + \sum_{i=1}^I w_i {\bf e}_{\bf u}^{(i)}, \delta p_h + \sum_{i=1}^I \gamma_i e_p^{(i)}\right).
\]

An alternative viewpoint on subspace correction algorithms comes from a linear-algebraic lens.  Here, we consider the discretized systems in~\eqref{eq:StokesSaddlePointSys} or~\eqref{eq:DiscNavStokesSaddlePointSys}, and define matrix restriction operators, $R^{(i)}$, that map from vectors in $\mathbb{R}^{n+m}$ to vectors in $\mathbb{R}^{n_i+m_i}$, where $n_i$ is the dimension of ${\bf V}_h^{(i)}$ and $m_i$ is the dimension of $W_h^{(i)}$.  Then, for the Stokes case, the equivalent updates to those in~\eqref{eq:UpdateStokes} are given by the solution of the restricted linear system
\begin{equation}\label{eq:UpdateStokes_LA}
R^{(i)} \mathcal{A} \left(R^{(i)}\right)^T {\bf x}^{(i)} = R^{(i)} \left({\bf b} - \mathcal{A}{\bf x}\right),
\end{equation}
and the updated discrete solution is given by
\[
{\bf x} \leftarrow {\bf x} + \sum_{i=1}^I \left(R^{(i)}\right)^T W^{(i)} {\bf x}^{(i)},
\]
where $W^{(i)}$ is the matrix representation of the weighting operators, $w_i$ and $\gamma_i$.  A similar update is given in the Navier-Stokes case by replacing matrix $\mathcal{A}$ and vectors ${\bf x}$ and ${\bf b}$ in~\eqref{eq:UpdateStokes_LA} by matrix $\tilde{\mathcal{A}}$ and vectors $\tilde{\bf x}$ and $\tilde{\bf b}$ from~\eqref{eq:DiscNavStokesSaddlePointSys}.

The descriptions above focus on applying the subspace correction algorithm directly to the linear systems discretized on the finest level.  However, in practice, the algorithm can be applied at any level of the multigrid hierarchy, and at any stage of the outer Newton iteration.  Thus, Algorithm~\ref{alg:ParSubCorrect} presents the algorithm for the Stokes equations from the linear-algebraic perspective with a generic right-hand side, given by ${\bf r}$, and a generic current approximation, ${\bf x}$, noting that these may, in fact, be approximations to the Newton updates or to coarse-grid quantities, depending on when the algorithm is called.  Since we focus on the use of FGMRES-accelerated relaxation, the right-hand side is generally expected to be an Arnoldi vector, while the current approximation is typically taken to be a zero vector; see~\cite[Remark 11.2]{AdlerSterckMacLachlanOlson}.  The same substitutions of $\mathcal{A} \rightarrow \tilde{\mathcal{A}}$, ${\bf x} \rightarrow \tilde{\bf x}$, and ${\bf b} \rightarrow \tilde{\bf b}$ are needed to express the algorithm applied to the Navier-Stokes case.

\begin{algorithm}[t]
\caption{ ({\bf Single sweep of additive subspace correction method for the Stokes equations})}
\label{alg:ParSubCorrect}
\algsetup{linenodelimiter=.}
{\bf Inputs}: 
\begin{itemize}
  \setlength\itemsep{-1pt}
\item matrix, $\mathcal{A}$
  \item subspace restriction operators, $R^{(i)}$
\item current approximate solution, ${\bf x}$
\item right-hand side, ${\bf r}$
\item weighting operators, $W^{(i)}$

\end{itemize}
{\bf Output}: Updated approximate solution, ${\bf x}$
\begin{algorithmic}[1]
\FOR {$i=1$ to $I$}
\STATE Find the solution, ${\bf x}^{(i)}$, of the local problem 
\[
R^{(i)} \mathcal{A} \left(R^{(i)}\right)^T {\bf x}^{(i)} = R^{(i)} \left({\bf r} - \mathcal{A}{\bf x}\right),
\]
\ENDFOR
\STATE $\displaystyle {\bf x} \leftarrow {\bf x} + \sum_{i=1}^I \left(R^{(i)}\right)^T W^{(i)} {\bf x}^{(i)}$
\RETURN ${\bf x}$
\end{algorithmic}
\end{algorithm}

Within monolithic multigrid solvers, it is common to design relaxation schemes where the subspace decomposition is inherited from the topology of the mesh.  Vanka's original decomposition for the 2D MAC-scheme finite-difference discretization used overlapping cell-wise patches, forming a patch for each cell in the mesh that included the cell-centred pressure DoF, as well as the four face-centred velocity DoFs~\cite{Vanka1}.  A similar construction was used for Crouzeix-Raviart elements in~\cite{NavierStokesVolkerTobiska}.  For discretizations with nodal pressure DoFs, such as the lowest-order Taylor-Hood case, John and Matthies proposed ``pressure-node oriented'' relaxation~\cite{JohnPatchesStokes}, where each patch includes a single pressure DoF and all algebraically connected velocity DoFs, where algebraic connection was defined based on nonzero entries in each row of $B$.  Topological variants of such approaches are implemented in PCPATCH~\cite{PCPATCH}, where decompositions can be formed based on the mesh topology.  For example, for the lowest-order Taylor-Hood elements, Vanka patches can be formed by taking each nodal pressure DoF and all velocity DoFs on the closures of the elements adjacent to the central node.  When no coincidental zeros appear in the matrix (i.e., for meshes lacking certain symmetries), this coincides with the pressure-node oriented approach from~\cite{JohnPatchesStokes}.  In~\cite{HighOrderASMStokes}, we proposed a composite Vanka preconditioner (explained in detail below) making use of such topological construction that extends the robustness of the pressure-node oriented Vanka relaxation to higher-order (generalized) Taylor-Hood elements.

An additional consideration explored in~\cite{HighOrderASMStokes} is the efficiency of two natural implementations of the subspace correction framework for higher-order finite-element spaces.  The algorithms implemented in PCPATCH~\cite{PCPATCH} make use of local assembly callbacks, to directly assemble the patch linear systems $R^{(i)} \mathcal{A} \left(R^{(i)}\right)^T$, without requiring the global assembly of $\mathcal{A}$.  This makes it possible to use subspace correction relaxation schemes for multigrid preconditioners within a matrix-free setting, albeit at the expense of independently assembling each patch matrix.  In contrast, the ASMPatchPC approach studied in~\cite{HighOrderASMStokes} requires the global assembly of $\mathcal{A}$, and forms the patches by extracting the local patch matrices from the global one.  The numerical experiments in~\cite{HighOrderASMStokes} suggest that this is a more efficient approach for higher-order discretizations of the Stokes equations with the relaxation scheme proposed there.  This can be attributed to the high cost of assembly over each element in the higher-order case, noting that the callback approach assembles the same elemental matrices several times (once for each patch in which each element appears), while the global assembly approach assembles over each element just once.  Since element-wise assembly cannot be avoided in patch-based relaxation schemes like these, minimizing the number of such assemblies seems naturally most efficient.  Here, we consider only the extraction approach implemented in ASMPatchPC, partly because of the conclusions from~\cite{HighOrderASMStokes} and partly because PCPATCH does not currently support the construction of local matrices for the subspace decomposition proposed next.

\subsection{Composite Vanka and Vanka-star relaxation schemes}\label{ssec:VankaStar}

\begin{figure}[t]
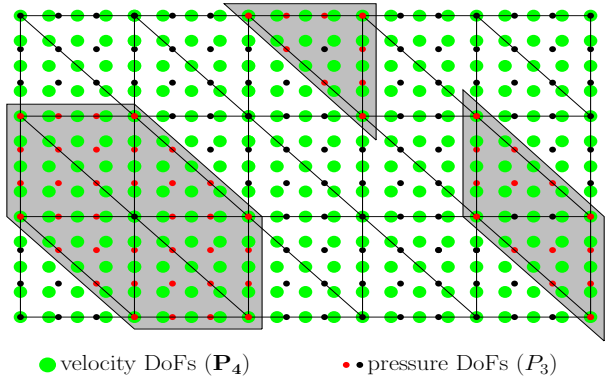
 
\begin{center}

\end{center}
\caption{Composite Vanka patch construction for triangular meshes. Two-dimensional triangular finite-element mesh for the $({\bf P_4},P_3)$ Taylor-Hood discretization with composite Vanka relaxation patches highlighted for vertex (left), edge (right), and element (top) patches. Large green circles denote velocity DoFs.  Red stars and small black circles show the pressure DoFs, where each patch contains all of the velocity DoFs in the block and the black circle pressure DoFs, but not the red star pressure DoFs within the patch.}
\label{fig:CompositeVankaPatchNavierStokes}
\end{figure}

In~\cite{HighOrderASMStokes}, we proposed a $p$-robust monolithic multigrid scheme for higher-order (generalized) Taylor-Hood discretizations of the Stokes equations in 2D, that made use of the topological generalization of the ``pressure-node oriented'' Vanka relaxation from John and Matthies~\cite{JohnPatchesStokes}.  In this method, we defined Vanka patches corresponding to each topological entity in the mesh.  In two dimensions, this leads to three types of patches covering all of the pressure DoFs, with patches for each of the nodes of the mesh (containing a single pressure DoF for a degree $p$ pressure space), each of the edges of the mesh (with $p-1$ pressure DoFs for each edge when $p>1$), and one for each of the elements of the mesh (with $(p+1)(p+2)/2 - 3p$ pressure DoFs when $p>2$).  For each patch, we take all velocity DoFs on the closure of the elements adjacent to the topological entity at the center of the patch.  This always results in all of the velocity DoFs on the closure of a single element for the element patches, and all of the velocity DoFs on the closure of two adjacent elements for the edge patches (except for boundary edges, which have only one adjacent element).  Vertex-based patches can, of course, vary in size, but there are six adjacent elements for the canonical triangularization of a regular mesh of nodes in 2D, as depicted in Figure~\ref{fig:CompositeVankaPatchNavierStokes}.

While the performance of monolithic multigrid with this relaxation gives clear robustness to $h-$ and $p-$refinement of regular grids, we note that the iteration counts observed in~\cite{HighOrderASMStokes} are somewhat higher than ideal for claiming ``textbook'' multigrid efficiency.  In this work, we step away from the algebraic and topological variants of ``pressure-node-oriented'' Vanka relaxation to propose a new patch construction that shows significant improvements in both iterations and time-to-solution over the composite relaxation scheme from~\cite{HighOrderASMStokes}.  Part of this improvement comes from a strong reduction in the number of patches, as we will go from a composite Vanka relaxation scheme, with one patch for each node, edge, and element in the (two-dimensional) mesh to a patch-based scheme that will have only one patch per vertex in the mesh.

In order to achieve effective multigrid relaxation with only vertex-based patches, it is clear that all pressure DoFs (as well as all velocity DoFs) must belong to at least one patch.  To our knowledge, we propose the first subspace decomposition whereby there is overlap between the pressure DoFs in the subdomains as well as with the velocity DoFs, breaking the ``pressure-node-oriented'' paradigm.  Quite simply, we propose to keep the velocity decomposition used for the vertex-centred patches in the composite Vanka scheme described above, but to enlarge the pressure space for each patch from a single pressure DoF to include all of the pressure DoFs on adjacent topological entities (but not their closures) to the central vertex.  In the topological patch description of PCPATCH~\cite{PCPATCH}, this is phrased as the DoFs in the \emph{star} of the vertex, so we call these patches \emph{Vanka-star} patches, to emphasize their essential nature as Vanka relaxation schemes, but based on vertex stars for the pressure, rather than just the single vertex-based pressure DoF.

\begin{figure}[t]
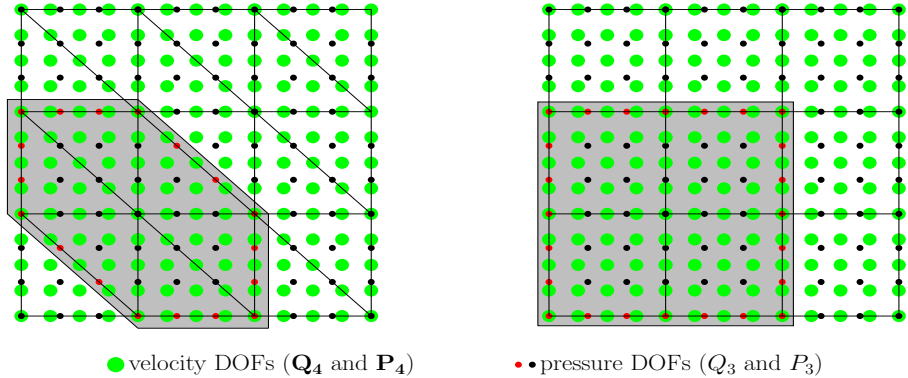
 
\begin{center}

\end{center}
\caption{Vanka-star patches for triangular and quadrilateral meshes.  {\bf Left}: Two-dimensional triangular finite-element mesh corresponding to $({\bf P_4},P_3)$ Taylor-Hood discretizations with vertex Vanka-star patch.
{\bf Right}: Quadrilateral finite-element mesh associated with $({\bf Q_4},Q_3)$ Taylor-Hood discretization with vertex Vanka-star patch. Large green circles are associated with velocity DoFs.  Red stars and small black circles show the pressure DoFs. Each Vanka-star patch contains all of the velocity DoFs in the highlighted block as well as the black circle pressure DoFs, but not the red star pressure DoFs within the patch. }
\label{fig:VankaStarPatchNavierStokes}
\end{figure}

Intuition into the construction of these patches comes from comparing the vertex-star patch shown for the triangular grid case in Figure \ref{fig:VankaStarPatchNavierStokes} to the composite patches shown in Figure \ref{fig:CompositeVankaPatchNavierStokes}.  Considering the edge and element patches for those entities adjacent to the central vertex in the vertex-based patch from the composite Vanka scheme, we can easily see that the velocity DoFs for each of these patches are a subset of those included in the vertex-based Vanka patch.  Thus, by adding the pressure DoFs for these entities to the pressure subspace in the subspace decomposition, we are maintaining the central idea of Vanka-style relaxation, that each patch contain some subset of the pressure DoFs along with all of the velocity DoFs that are connected via the rows of $B$.  The given construction also naturally extends to three-dimensional meshes, where we can easily define both the star of a vertex for the pressure DoFs and the closure of the star for velocity DoFs to give the resulting construction.

\section{2D Numerical Experiments}\label{sec:2Dnumerics}

In this section, we measure the efficiency of the proposed preconditioner for solution of the discretized two-dimensional Stokes and Navier-Stokes problems using test problems on the unit-square domain.  All numerical experiments were performed in Firedrake~\cite{Firedrake}, making use of the tight integration with PETSc~\cite{petsc} for the linear and nonlinear solvers~\cite{kirby2018solver}.  All experiments in this section were performed using 8 cores of a server with dual 8-core Intel Xeon 1.7 GHz CPUs and 384 GB of RAM.

For all 2D experiments, we start with a $5\times 5$ coarsest grid of the unit square.  For quadrilateral meshes, we simply refine this uniformly to create multigrid hierarchies with $\ell$ levels, so that $\ell = 1$ corresponds to a $10\times 10$ mesh, and $\ell = 4$ corresponds to an $80\times 80$ mesh.  For triangular grids, we first cut each square element of the $5\times 5$ mesh into two triangles, dividing from top-left to bottom-right, then use uniform refinement on the resulting triangulation, leading to structured triangulations on all levels of the multigrid hierarchy, with $200$ elements for $\ell = 1$, and $12,800$ elements for $\ell = 4$.  The total numbers of degrees of freedom that arise in such discretizations grows rapidly with $k$, as shown in Table~\ref{table:nDoFs2DStokes2DNavierStokes}.
\begin{table}[ht]
\begin{center}
\captionsetup{justification=centering,margin=0.2cm}
\caption{{\small Number of DoFs for the finest-grid discretization for 2D (Navier-)Stokes problems with $({\bf P_k},P_{k-1})$ or $({\bf Q_k},Q_{k-1})$ discretizations and $\ell$ levels of refinement.}}
\label{table:nDoFs2DStokes2DNavierStokes}
{\small
\begin{tabular}{|c||c|c|c|c|c|c|}
        \hline
	\diagbox[]{$\ell$}{$k$} & 2 & 3 & 4 & 5 & 6 & 7\\
        \hline
        1 & \num{1003}    & \num{2363}  & \num{4323} & \num{6883} &  \num{10043} & \num{13803} \\  \hline
        2 & \num{3803}    & \num{9123}  & \num{16843} & \num{26963} & \num{39483}  & \num{54403} \\  \hline
        3 & \num{14803}  & \num{35843}  & \num{66483} & \num{106723} & \num{156563}  & \num{216003} \\  \hline
        4 & \num{58403}  & \num{142083}  & \num{264163} & \num{424643} & \num{623523} & \num{860803} \\  \hline
    \end{tabular}
}
\end{center}
\end{table}
In both cases, we use the canonical finite-element interpolation operators to interpolate from level $\ell-1$ to level $\ell$ (depending on the polynomial order, $k$, of the discretization).  For all cases, we use a zero initial guess to the solution, and perform up to 50 iterations of the nonlinear solver (for the Navier-Stokes case) and up to 100 iterations of FGMRES with no restart, marking problems for which convergence was not obtained as solver failures.  For this, we require reduction of the linear or nonlinear residual norm below a relative tolerance of $10^{-10}$.  For the Navier-Stokes equations, we use the Eisenstat-Walker method of setting linear solver tolerances for each linearization~\cite{EisenStatWalker}.

\begin{figure}[H]
\graphicspath{{./Figures/}}
\begin{center}
\includegraphics[width=0.7\textwidth]{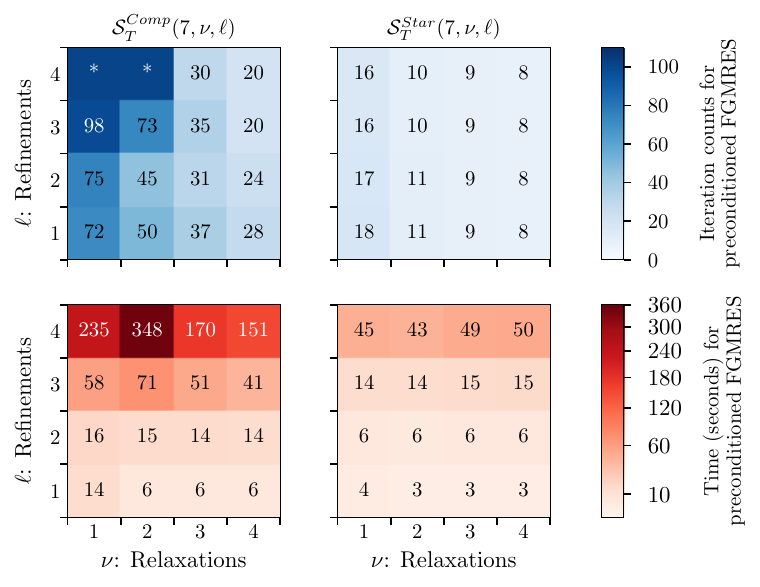}
\end{center}
\caption{Comparison of iteration counts ({\bf top}) and iteration time ({\bf bottom}) for the $({\bf P_7},P_{6})$ discretization of the Stokes equations, using monolithic-multigrid preconditioned FGMRES with composite Vanka relaxation ({\bf left}) and Vanka-star relaxation ({\bf right}). Results denoted by $\ast$ indicate a failure to converge in 100 FGMRES iterations, with time reported for those 100 iterations.}
\label{fig:P7P6_ASMPatchPC_CopmosVanka_vs_VankaStar_It_Time_FGMRES}
\end{figure}

\subsection{Stokes Equations}

Here, we consider a simple analytical solution for the velocity and pressure given by
\begin{eqnarray*}\label{eq:ArtificialStokes}
{\bf u}(x,y) = \left( \begin{array}{c}
4x^2y \\
-4xy^2
\end{array}\right), & \quad & p(x,y) = xy - \frac{1}{4},
\end{eqnarray*}
substituting these into~\eqref{eq:Stokes1} to compute the forcing function, ${\bf f}$, and corresponding Dirichlet boundary data.

We use $\mathcal{S}^{Comp}_{T}(k, \nu, \ell)$ ($\mathcal{S}^{Star}_{T}(k, \nu, \ell)$) to denote the composite Vanka (and Vanka-star) based solvers using monolithic $V(\nu,\nu)$ cycles as preconditioners for FGMRES(100) with $1\leq \ell \leq 4$ levels of refinement on triangular meshes using the $({\bf P_k},P_{k-1})$ discretization for $k\ge 2$.  For solvers on quadrilateral meshes, using the $({\bf Q_k},Q_{k-1})$ discretization for $k\ge 2$, we replace the subscript $T$ by $Q$.  In the experiments that follow, we will focus on robustness of the solvers to variations in $k$ and $\ell$, and the impact of $\nu$ on the cost of solution using the two different relaxation schemes.

Figure~\ref{fig:P7P6_ASMPatchPC_CopmosVanka_vs_VankaStar_It_Time_FGMRES} presents iteration counts and time-to-solution for the case of $k=7$, focusing on the comparison between performance using the composite Vanka relaxation scheme (at left) and the Vanka-star scheme (at right).  Immediately apparent is the significant reduction in the number of iterations needed for convergence when using Vanka-Star relaxation in comparison to composite Vanka relaxation.  In the top-right graph, we see greatly improved scalability in the iteration counts with number of levels of refinement, $\ell$, particularly when using $\nu = 1$ or $2$ relaxation sweeps per V-cycle.  This improved performance is also reflected in the time-to-solution.  Looking at the finest grid, $\ell = 4$, we see that the fastest time-to-solution when using Vanka-star relaxation is 43 seconds (with $\nu = 2$), over 3 times faster than the best time-to-solution using composite Vanka (151 seconds with $\nu = 4$).  A similar speedup is observed for the next-finest grid, $\ell = 3$, where fewer convergence issues are observed with composite Vanka relaxation.

\begin{figure}[t]
\graphicspath{{./Figures/}}
\begin{center}
\includegraphics{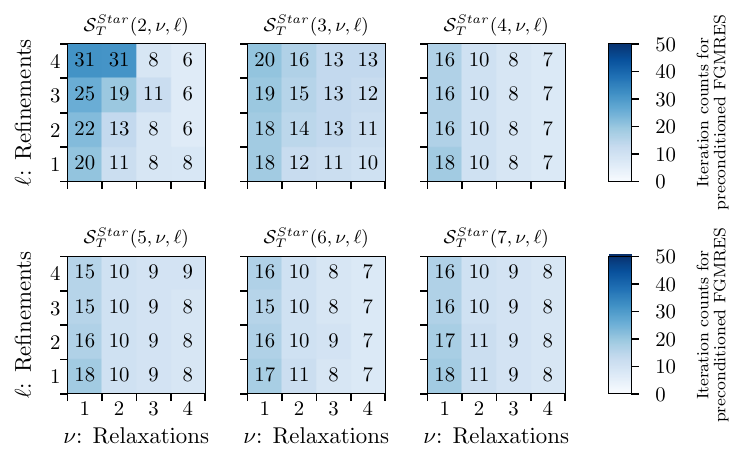}
\end{center}
\caption{Iteration counts for the $({\bf P_k},P_{k-1})$ discretization of the Stokes equations using Vanka-star relaxation, as we increase the discretization order, $k$, from 2 to 7.}
\label{fig:All_PP_ASMPatchPC_VankaStar_ItFGMRES}
\end{figure}

In Figure \ref{fig:All_PP_ASMPatchPC_VankaStar_ItFGMRES}, we consider the robustness of the iteration counts to solution as we vary the order of the discretization, $k$.  For lower orders ($k=2$ and $3$), we see somewhat higher iterations for $\nu=1$ and $2$ but, in all cases, iteration counts are reasonable, and show expected behaviour as we increase the number of relaxation sweeps per V-cycle.  For $k\geq 4$, we see very little variation in the iteration counts with order.  

\begin{figure}[H]
\graphicspath{{./Figures/}}
\begin{center}
\includegraphics{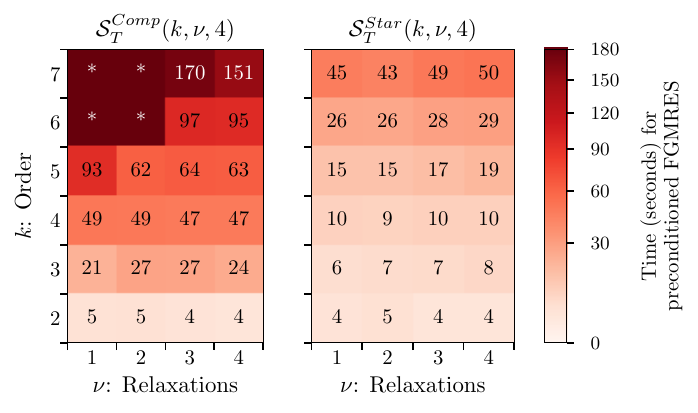}
\end{center}
\caption{Comparison of iteration times for solvers using composite Vanka relaxation ({\bf left}) and Vanka-star relaxation ({\bf right}) as we vary the polynomial order, $k$, and number of relaxation sweeps used, $\nu$, on triangular meshes.  Results denoted by $\ast$ indicate a failure to converge in 100 FGMRES iterations.}
\label{fig:All_Timings_PP_nref4_ASMPatchPC_ComposVanka_VankaStar}
\end{figure}

Figure \ref{fig:All_Timings_PP_nref4_ASMPatchPC_ComposVanka_VankaStar} focuses on CPU timings for the finest grid, $\ell = 4$, as we vary polynomial order, $k$, and number of relaxation sweeps per $V$-cycle, $\nu$.  Considering the Vanka-star relaxation, with data at right, we see very little impact on time-to-solution as we vary $\nu$, with best time generally for $\nu = 1$ or $2$, but only a few seconds of variation.  Moreover, the cost is seen to increase slowly with $k$.  Comparing the best time-to-solution for $k=3$ and $k=6$, we see an increase by a factor of just over 4, suggesting better-than-expected time-to-solution scaling like $\mathcal{O}(k^2)$ (noting that the patch matrices grow in size like $\mathcal{O}(k^2)$, so standard Gaussian elimination would scale like $\mathcal{O}(k^6)$ if this were the dominant cost of the iteration).  Comparing performance between the solver using composite Vanka relaxation (at left) and Vanka-star (at right), we see substantial improvements with Vanka-star relaxation, showing speedups of up to four times for $k>2$ (noting that at $k=2$, the two approaches coincide).

\begin{figure}[t]
\graphicspath{{./Figures/}}
\begin{center}
\includegraphics[width=0.7\textwidth]{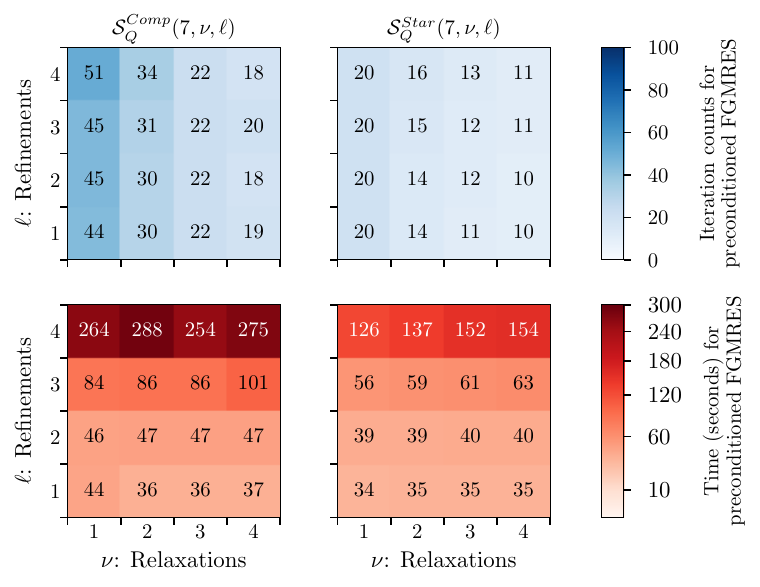}
\end{center}
\caption{Comparison of iteration counts ({\bf top}) and iteration time ({\bf bottom}) for the $({\bf Q_7},Q_{6})$ discretization of the Stokes equations, using monolithic-multigrid preconditioned FGMRES with composite Vanka relaxation ({\bf left}) and Vanka-star relaxation ({\bf right}).}
\label{fig:Q7Q6_ASMPatchPC_CopmosVanka_vs_VankaStar_It_Time_FGMRES}
\end{figure}

For the quadrilateral case, Figure~\ref{fig:Q7Q6_ASMPatchPC_CopmosVanka_vs_VankaStar_It_Time_FGMRES} shows the same comparison as in Figure~\ref{fig:P7P6_ASMPatchPC_CopmosVanka_vs_VankaStar_It_Time_FGMRES} between solvers using composite Vanka and those using Vanka-star relaxation for the  $({\bf Q_7},Q_{6})$ discretization.  While the performance of the solver using composite Vanka relaxation is much better in this case than for triangular grids (as noted in~\cite{HighOrderASMStokes}), the comparison is, nonetheless, quite similar.  Again, the use of Vanka-star relaxation leads to significant improvements in iterations to solution, resulting in notable speedups over composite Vanka.  Now, since the composite Vanka results are not quite so bad, the speedups are more modest, but still reach a factor of two on the finest grid, with $\ell = 4$.  

Comparing time-to-solution between the lower-right data in Figure~\ref{fig:Q7Q6_ASMPatchPC_CopmosVanka_vs_VankaStar_It_Time_FGMRES} with that in Figure~\ref{fig:P7P6_ASMPatchPC_CopmosVanka_vs_VankaStar_It_Time_FGMRES}, we see that the triangular grid discretization has a much faster solve time, due to a combination of slightly lower iteration counts and faster time-per-iteration.  The faster time-per-iteration is to be expected, because we have the same number of patches in these two cases, but those for the triangular-grid discretization have fewer DoFs per patch, due to their construction.  We omit the analogue of Figure~\ref{fig:All_PP_ASMPatchPC_VankaStar_ItFGMRES} for quadrilateral case, because it is very similar in form to that figure, albeit with iteration counts for $k\geq 4$ that are close to those shown in the top-right data from Figure~\ref{fig:Q7Q6_ASMPatchPC_CopmosVanka_vs_VankaStar_It_Time_FGMRES}.

\begin{figure}[H]
\graphicspath{{./Figures/}}
\begin{center}
\includegraphics{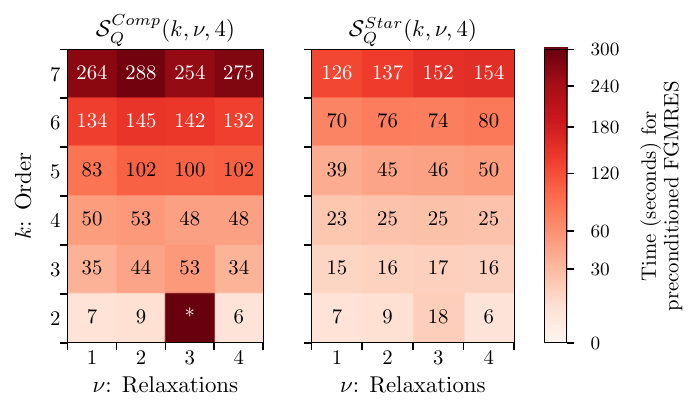}
\end{center}
\caption{Comparison of iteration times for solvers using composite Vanka relaxation ({\bf left}) and Vanka-star relaxation ({\bf right}) as we vary the polynomial order, $k$, and number of relaxation sweeps used, $\nu$ on quadrilateral meshes.  Results denoted by $\ast$ indicate a failure to converge in 100 FGMRES iterations.}
\label{fig:All_Timings_QQ_nref4_ASMPatchPC_ComposVanka_VankaStar}
\end{figure}

Figure~\ref{fig:All_Timings_QQ_nref4_ASMPatchPC_ComposVanka_VankaStar} provides the same comparison in time-to-solution as we vary polynomial order, $k$, and number of relaxation sweeps per $V$-cycle, $\nu$, as in Figure~\ref{fig:All_Timings_PP_nref4_ASMPatchPC_ComposVanka_VankaStar}, but for quadrilateral meshes instead of triangular ones.  Aside from an isolated solver failure at lowest order (that we did not investigate further), the conclusions from these results are very similar to those from the triangular case.  At high orders, the best solvers using Vanka-star relaxation are consistently about twice as fast as those using composite Vanka relaxation.  The best results for Vanka-star are consistently obtained using $\nu=1$ relaxation sweeps per $V$-cycle.  Comparing times between the triangular and quadrilateral grid cases, we again see that solvers at the same order on triangular grids are faster than those on quadrilateral grids, due to lower numbers of iterations to convergence and lower costs per iteration.

\subsection{Navier-Stokes equations}

Here, we make use of the standard 2D Navier-Stokes lid-driven cavity test problem, setting the forcing term, $\mathbf{f} = \mathbf{0}$, with homogeneous Dirichlet boundary conditions on three faces of the unit-square domain, and $\mathbf{u} = (1,0)^T$ imposed on the top face.  We focus on the Newton-Krylov-Multigrid methodology here, using the Eisenstat-Walker criteria~\cite{EisenStatWalker}, with default parameters in PETSc, to choose the linear stopping criteria for the monolithic-multigrid-preconditioned FGMRES iterations used to solve each Newton linearization.  Since performance now depends on both the solver parameters and the Reynolds number, $Re$, of the flow, we now denote the solver configuration by $\mathcal{NS}^{Star}_{T}(Re,k, \nu, \ell)$, possibly with superscript $Comp$ to denote using composite Vanka relaxation in place of Vanka-star relaxation, or with subscript $Q$ to denote quadrilateral meshes in place of triangular ones.

\begin{figure}[H]
\graphicspath{{./Figures/}}
\begin{center}
\includegraphics[width=0.8\textwidth]{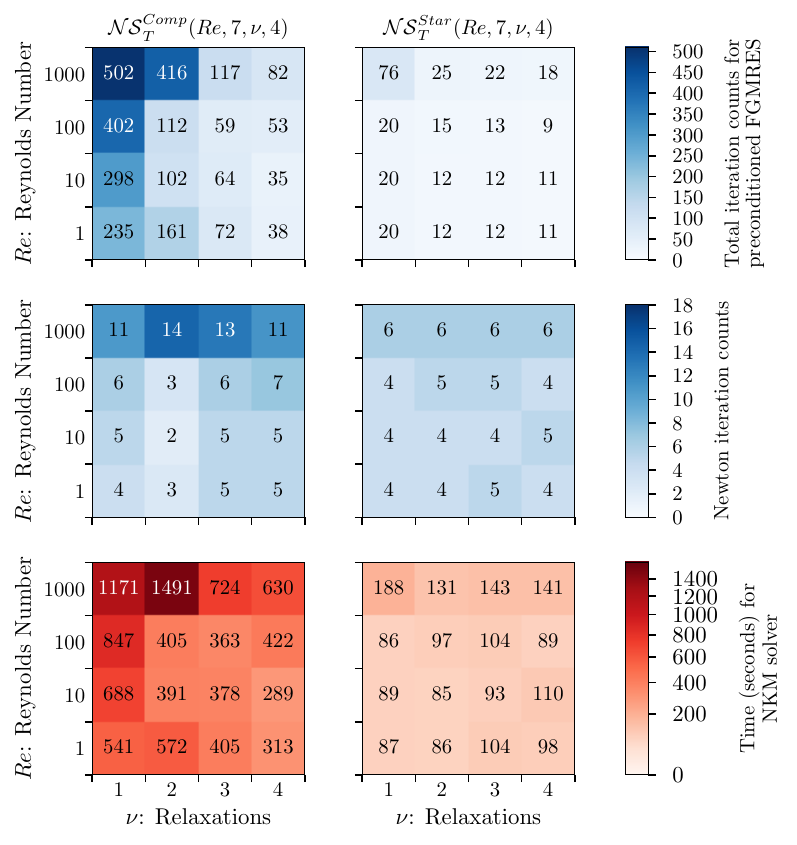}
\end{center}
\caption{Comparison of total FGMRES iteration counts ({\bf top}), Newton iteration counts ({\bf middle}), and total solve time ({\bf bottom row}) for the $({\bf P_7},P_{6})$ discretization of the Navier-Stokes equations on triangular grids, using Newton-Krylov-Multigrid solvers with composite Vanka relaxation ({\bf left}) and Vanka-star relaxation ({\bf right}).}
\label{fig:NST2_k7nref4_ASM_CompStarVanka_LIt_NLIt_Time_FGMRES}
\end{figure}

Figure \ref{fig:NST2_k7nref4_ASM_CompStarVanka_LIt_NLIt_Time_FGMRES} presents a comparison between the solvers using composite Vanka (at left) and Vanka-star (at right) relaxation for the $({\bf P_7},P_{6})$ discretization at grid resolution $\ell = 4$, as we vary $Re$ over 3 orders of magnitude, using $1 \leq \nu \leq 4$ relaxation sweeps within the V-cycle.  The top row of this figure presents the total number of linear iterations over all nonlinear iterations, while the total number of nonlinear iterations is presented in the middle row of figure, and the total CPU time-to-solution is presented in the bottom row.  We expect the problems at $Re = 1$ to be relatively easy, while those at $Re = 1000$ should be quite challenging at this grid resolution, even for a higher-order discretization.  That this is true is immediately apparent in all rows, with generally increasing linear and nonlinear iteration counts and solution times as $Re$ increases.  Comparing linear iteration counts, we see substantially more required for the solvers using composite Vanka relaxation in all cases than for those using Vanka-star relaxation.  This is partly due to the also generally increased 
\begin{figure}[H]
\graphicspath{{./Figures/}}
\begin{center}
\includegraphics[width=0.7\textwidth]{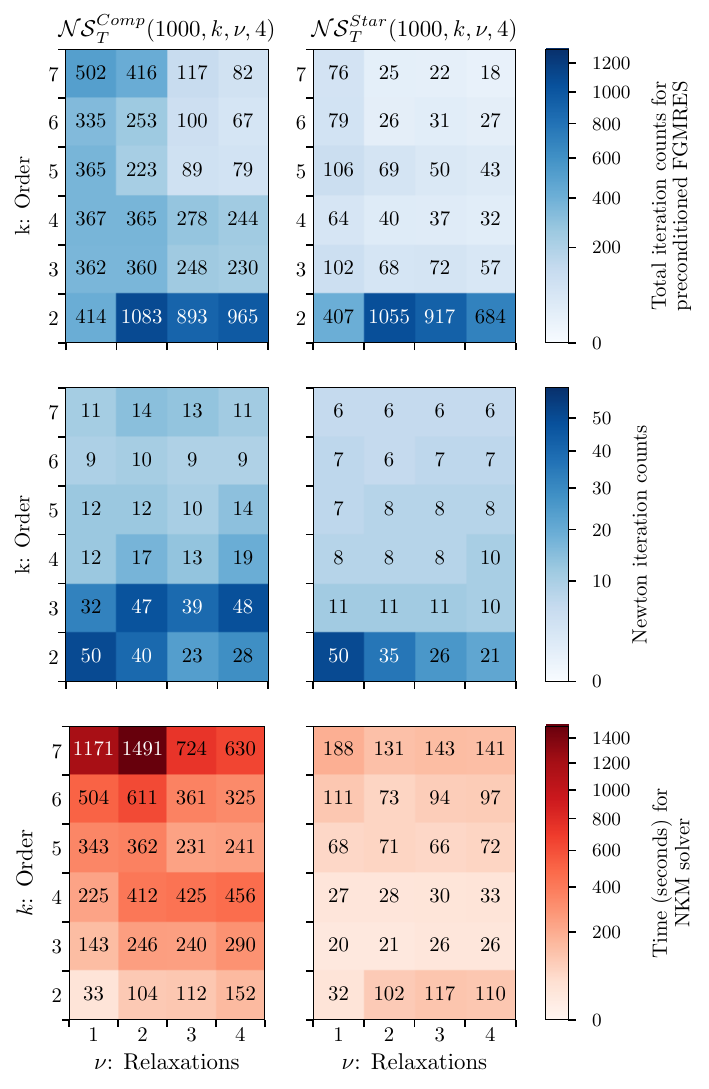}
\end{center}
\caption{Comparison of total FGMRES iteration counts ({\bf top}), Newton iteration counts ({\bf middle}), and total solve time ({\bf bottom row}) for the $({\bf P_{k}},P_{k-1})$ discretization of the Navier-Stokes equations at $Re = 1000$ on triangular grids, using Newton-Krylov-Multigrid solvers with composite Vanka relaxation ({\bf left}) and Vanka-star relaxation ({\bf right}).}
\label{fig:NST2_Re1000nref4_ASM_CompStarVanka_LIt_NLIt_Time_FGMRES}
\end{figure}
\noindent number of Newton iterations needed for the solver using composite Vanka relaxation, noting that this is likely impacted by the use of the Eisenstat-Walker stopping criteria, which are impacted by the quality of the preconditioner for the linearizations.  Overall, it is quite apparent that the linear and nonlinear iteration counts for the solver using Vanka-star relaxation are much more consistent as we vary $Re$, with only one outlier in the data, where the case of $V(1,1)$ cycles with $Re = 1000$ requiring substantially more linear iterations than any other case.  
This steadiness in iterations is also reflected in the time-to-solution data, where we see much lower solution times and much less variation in solution times for the solver using Vanka-star relaxation than for that using composite Vanka.  
Even with $Re = 1$, we see almost a $4\times$ speedup with Vanka-star relaxation (comparing best solution times from that row of the figure), which increases to almost $5\times$ speedup at $Re = 1000$.

We next fix $Re = 1000$ and $\ell = 4$, and study the impact of the discretization order, $k$, on the performance of the solvers. Figure \ref{fig:NST2_Re1000nref4_ASM_CompStarVanka_LIt_NLIt_Time_FGMRES} reports the same data as above for this case, again considering the impact of the number of relaxation sweeps within the monolithic V-cycle preconditioner.  Here, with fixed (large) $Re$, we expect the most difficult solves to be at low order, since the $({\bf P_k},P_{k-1})$ discretization does a poor job of resolving the flow at this Reynolds number with small $k$.  This is, indeed, reflected in the solver statistics reported, where we note that we allowed a maximum of 50 nonlinear iterations in these results, so the results for both solvers with $k=2$ and $\nu = 1$ should be interpreted as solver failures.  Within this data, we note that using Vanka-star relaxation leads to both generally lower linear and nonlinear iteration counts, with much more graceful failure as $k$ decreases.  While no solver can be said to be performing well at $k=2$, the data for Vanka-star relaxation at $k=3$ or $4$ is not substantially different than at higher orders.  All of this is, once again, reflected in the time-to-solution data, where we see $4\times$  or better speedup for $k>5$, while we see bigger speedups, by factors of seven or eight, at $k=3$ or $4$.

\begin{figure}[t]
\graphicspath{{./Figures/}}
\begin{center}
\includegraphics{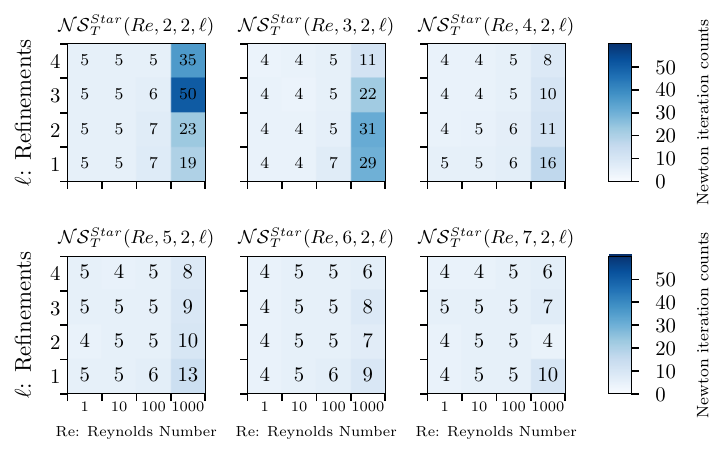}
\end{center}
\caption{Newton iteration counts for solving the $({\bf P_{k}},P_{k-1})$ discretization of the Navier-Stokes equations using Vanka-star relaxation with $\nu=2$ sweeps of relaxation, as a function of polynomial order, $k$, refinement level, $\ell$, and Reynolds number, $Re$.}
\label{fig:All_NST2_nu2_ASM_VankaStar_NLIt_Newton}
\end{figure}

\begin{figure}[H]
\graphicspath{{./Figures/}}
\begin{center}
\includegraphics{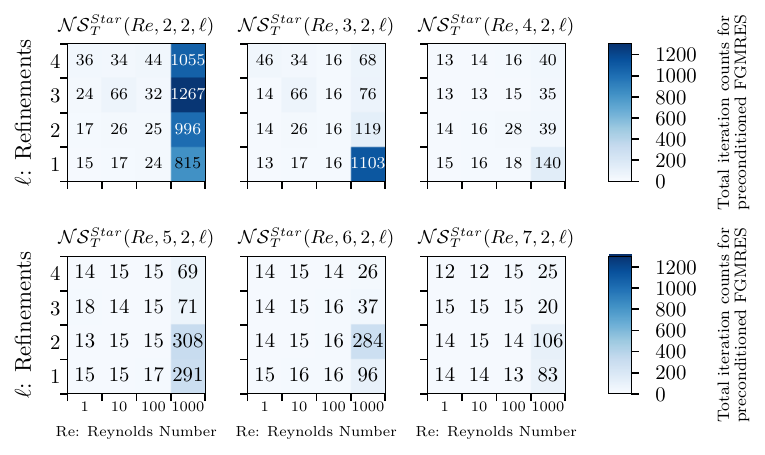}
\end{center}
\caption{Total FGMRES iteration counts for solving the $({\bf P_{k}},P_{k-1})$ discretization of the Navier-Stokes equations using Vanka-star relaxation with $\nu=2$ sweeps of relaxation, as a function of polynomial order, $k$, refinement level, $\ell$, and Reynolds number, $Re$.}
\label{fig:All_NST2_nu2_ASM_VankaStar_LIt_FGMRES}
\end{figure}

Figures~\ref{fig:All_NST2_nu2_ASM_VankaStar_NLIt_Newton} and~\ref{fig:All_NST2_nu2_ASM_VankaStar_LIt_FGMRES} present nonlinear and linear iteration counts, respectively, to solution with $\nu = 2$ as we change the discretization order, $k$, level of refinement, $\ell$, and Reynolds number $Re$.  Very little surprising occurs in the nonlinear iteration counts reported in Figure~\ref{fig:All_NST2_nu2_ASM_VankaStar_NLIt_Newton} except, perhaps, in the reasonable performance of the solvers at low levels of grid refinement when the Reynolds number is $100$, although this is probably indicative of convergence to a poor-quality solution, since the $10\times 10$ grid at $\ell =1$ cannot possibly resolve the flow with $Re = 100$ when $k$ is small.  Overall, aside from large Newton iteration counts at $Re = 1000$ when $k = 2$ or $3$, the nonlinear solver converges quite reasonably.  Performance of the linear solver is similar, although we see large number of linear iterations persisting for $\ell = 1$ or $2$ with $Re = 1000$ even for large $k$, which is to be expected since, again, the flow at this Reynolds number is underresolved on these grids.  Analogous figures for the solvers using composite Vanka and $\nu = 4$ (not shown here) show similar results but with consistently higher iteration counts, despite the added relaxation for each V-cycle.  While Newton iteration counts for this solver are comparable to those using Vanka-star relaxation for small $k$ and $Re \leq 100$, we saw consistent growth in nonlinear iteration counts for $Re = 1000$ at all discretization orders and at $Re = 100$ \noindent  on coarser grids. Linear iteration counts show greater disparity (as in the figures above), with typical total linear iteration counts in the range of 20-60, at discretization orders $k\ge 4$ and $Re\le 100$, for the solver using composite Vanka relaxation, in comparison to the 10-20 iterations for most similar cases in Figure~\ref{fig:All_NST2_nu2_ASM_VankaStar_LIt_FGMRES}.

\begin{figure}[H]
\graphicspath{{./Figures/}}
\begin{center}
\includegraphics[width=0.7\textwidth]{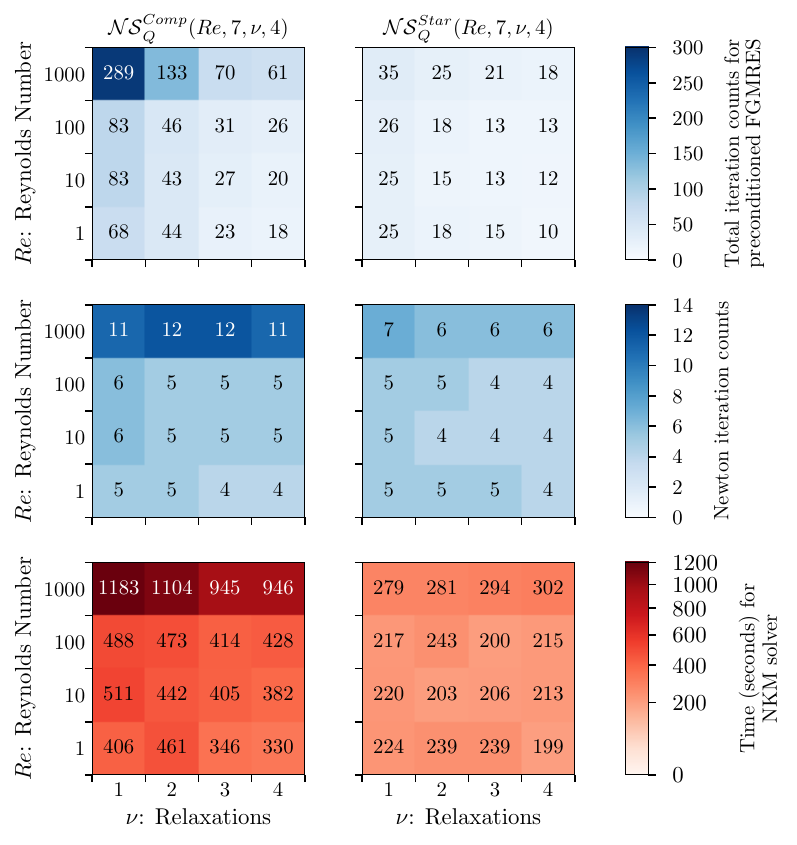}
\end{center}
\caption{Comparison of total FGMRES iteration counts ({\bf top}), Newton iteration counts ({\bf middle}), and total solve time ({\bf bottom row}) for the $({\bf Q_7},Q_{6})$ discretization of the Navier-Stokes equations on triangular grids, using Newton-Krylov-Multigrid solvers with composite Vanka relaxation ({\bf left}) and Vanka-star relaxation ({\bf right}).}
\label{fig:NSQ2_k7nref4_ASM_CompStarVanka_LIt_NLIt_Time_FGMRES}
\end{figure}

To complete this section, we demonstrate similar performance for the solvers using the two relaxation schemes for the discretized lid-driven cavity problem on quadrilateral meshes.  Figure~\ref{fig:NSQ2_k7nref4_ASM_CompStarVanka_LIt_NLIt_Time_FGMRES} repeats the experiment from Figure~\ref{fig:NST2_k7nref4_ASM_CompStarVanka_LIt_NLIt_Time_FGMRES}, just for the $({\bf Q_{7}},Q_{6})$ discretization in place of the $({\bf P_{7}},P_{6})$ discretization.  Here, we see somewhat improved performance of the solver using composite Vanka in comparison to the triangular grid case, albeit still with degradation at $Re = 1000$.  Comparing time-to-solution data from the bottom block of Figure \ref{fig:NSQ2_k7nref4_ASM_CompStarVanka_LIt_NLIt_Time_FGMRES}, we see that using the solver with Vanka-star relaxation leads to speedups by factors of 2 or 3 for higher Reynolds numbers, over the solver with composite Vanka, with somewhat more modest speedups at $Re = 1$.

\begin{figure}[H]
\graphicspath{{./Figures/}}
\begin{center}
\includegraphics[width=0.7\textwidth]{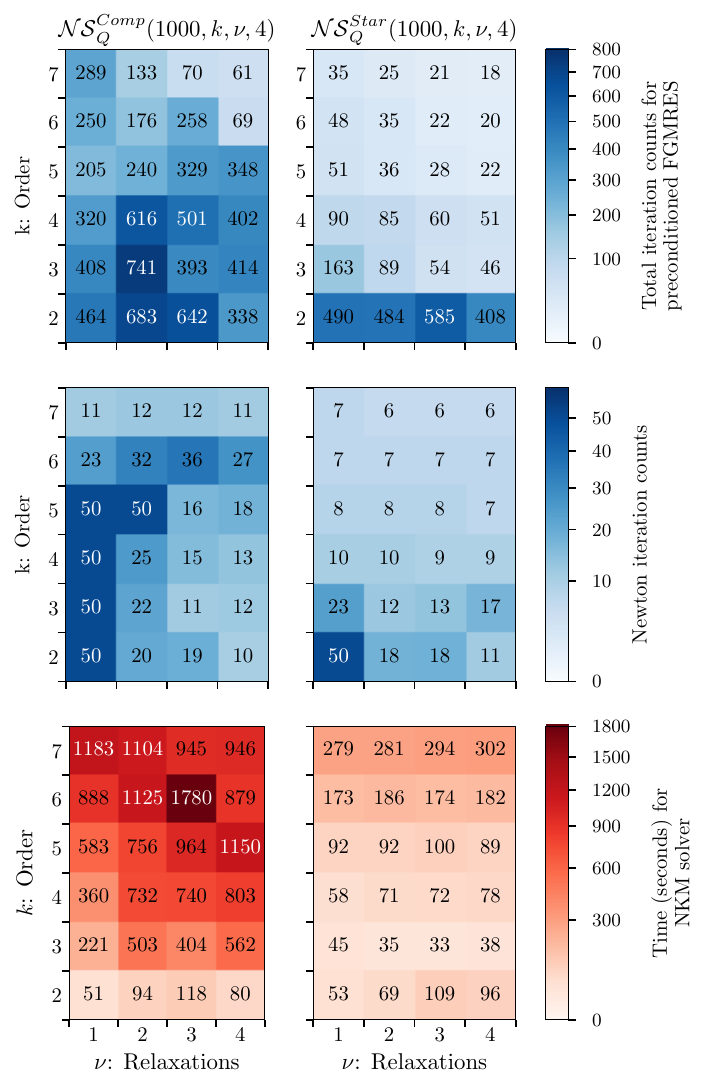}
\end{center}
\caption{Comparison of total FGMRES iteration counts ({\bf top}), Newton iteration counts ({\bf middle}), and total solve time ({\bf bottom row}) for the $({\bf Q_{k}},Q_{k-1})$ discretization of the Navier-Stokes equations at $Re = 1000$ on triangular grids, using Newton-Krylov-Multigrid solvers with composite Vanka relaxation ({\bf left}) and Vanka-star relaxation ({\bf right}).}
\label{fig:NSQ2_Re1000nref4_ASM_CompStarVanka_LIt_NLIt_Time_FGMRES}
\end{figure}

Figure~\ref{fig:NSQ2_Re1000nref4_ASM_CompStarVanka_LIt_NLIt_Time_FGMRES} presents analogous results to Figure~\ref{fig:NST2_Re1000nref4_ASM_CompStarVanka_LIt_NLIt_Time_FGMRES}, now for the $({\bf Q_k},Q_{k-1})$ discretization at $Re = 1000$.
Overall, we see similar behaviour for the solver using Vanka-star relaxation as we did above, and a similar comparison between the solvers on triangular and quadrilateral grids.  Here, we note worsening performance for the solver using composite Vanka relaxation (with many solver failures where the number of nonlinear iterations is reported as $50$), leading to similar increases in the number of linear iterations needed as well.  In comparison, the iteration counts for the solver with Vanka-star relaxation are quite steady beyond $k=4$.  We see similar improvements in time-to-solution for the solver using Vanka-star relaxation (again noting that some of the low times reported for the solver using composite Vanka relaxation are for runs that failed), with speedups of more than $7\times$ for low orders $k = 3, 4$, and speedups of $3$ to $6\times$ for high orders $k= 5, 6$.

\section{3D Numerical Experiments}\label{sec:3Dnumerics}
In this section, we demonstrate the effectiveness of the preconditioned FGMRES(30) solver for solution of the discretized three-dimensional (Navier-)Stokes problem. 
We use the proposed monolithic multigrid V-cycle preconditioner on tetrahedral finite-element meshes with the $({\bf P_{k}},P_{k-1})$ discretization on every level of the multigrid hierarchy. The relaxation at each level is three sweeps of Vanka-star preconditioned FGMRES, as both pre- and post-relaxation in all tests.  

\noindent We use direct solvers for the linear systems that arise from both the patch solves and the coarsest grid of the V-cycle.  
We use a test problem with manufactured solution for the Stokes case, and the lid-driven-cavity benchmark model problem 
\noindent as the Navier-Stokes test case, both defined on the unit-cube domain.  As in the 2D case, all numerical experiments were performed in Firedrake using PETSc within the linear and nonlinear solvers. For the lid-driven-cavity model, the inexact Newton linearization is coupled with the Eisenstat-Walker method for determining the solver tolerances.



These tests were run on the Digital Research Alliance of Canada cluster, Nibi, using a single node with two 2.4GHz Intel 6972P processors with 96 cores each and 6TB of RAM.
To form the coarsest mesh within the multigrid hierarchy for these problems, we first divide the unit cube into $4\times4\times4$ subcubes, each of which is then cut into six tetrahedra in the usual way.  These tetrahedral meshes are then refined $\ell$ times for $\ell = 2$ or $3$. 
For experiments with monolithic multigrid preconditioning with 2 and 3 levels of geometric refinement,
we use 6 and 48 cores, respectively, but must increase the memory allocated per core as we increase $k$, due to the increased size and density of the resulting linear and non-linear systems of equations. For the Stokes case, we experimentally increased the memory allocation for each core from 4 GB to 12 GB to 32 GB to 64 GB as we increased the discretization order within the multigrid hierarchy from $k=3$ to $6$. 
For the Navier-Stokes case, we experimentally found that memory allocation for each core must be increased from 4 GB to 15 GB to 50 GB to 120 GB as we increase the discretization order from $k=3$ to $k=4$ to $k=5$ to $k=6$.   Measured memory usage statistics are presented in Tables~\ref{table:PrecondSolverVankaStar3DStokes} and \ref{table:NKMVankaStar3DNavierStokes} below.  Table~\ref{table:nDoFs3DStokes3DNavierStokes} shows how the number of discrete DoFs grow with $k$ and $\ell$.

\begin{table}[t]
\begin{center}
\captionsetup{justification=centering,margin=0.2cm}
\caption{{\small Number of DoFs for the finest-grid discretization for 3D (Navier-)Stokes problems with $\ell$ levels of refinement from the $4\times 4\times 4$ base mesh and $({\bf P_k},P_{k-1})$ discretizations.}}
\label{table:nDoFs3DStokes3DNavierStokes}
{\small
\begin{tabular}{|c||c|c|c|c|}
\hline    
	\diagbox[]{$\ell$}{$k$} & 3 & 4 & 5 & 6 \\
        \hline
        2 &  \num{388884}    &  \num{941524}   & \num{1868948}   & \num{3269460} \\  \hline
        3 &  \num{3012644}  &  \num{7352740}  & \num{14666532} & \num{25740452} \\  \hline
    \end{tabular}
 }
\end{center}
\end{table}

\subsection{Stokes Equations}
Table \ref{table:PrecondSolverVankaStar3DStokes} shows numerical data associated with the solver for 3D Stokes equation with manufactured solution
\begin{eqnarray}\label{eq:3DManuFactSol}
{\bf u}(x,y,z) = \left( \begin{array}{c}2\sin(\pi x)\cos(\pi y)\cos(\pi z) \\
-\cos(\pi x)\sin(\pi y)\cos(\pi z) \\
-\cos(\pi x)\cos(\pi y)\sin(\pi z) \end{array}\right), &  &  p(x,y,z) = 0,
\end{eqnarray}
and Dirichlet boundary conditions imposed on ${\bf u}$ on all of boundary of the unit cube to match the manufactured solution.
We first hand-tune the relative stopping tolerance for the outer FGMRES solver, in order to ensure that we achieve the expected convergence to the level of finite-element discretization accuracy as we vary $k$.  For each value of $k$, we report these tolerances in the column labelled ``rtol'', and the measured ${\bf H}^1$-norm error in the velocity, $\|{\bf u}-{\bf u}_h\|_1$, and $L^2$-norm error in the pressure, $\|p-p_h\|_0$, for $\ell = 2$ and $3$ levels of refinement. In these columns, we see velocity and pressure error reductions between results for $\ell = 2$ and $\ell = 3$ by factors of about 8 for $k=3$ through about 64 for $k=6$, matching the expected finite-element convergence. We note that the required relative residual stopping tolerance to achieve this at $k=6$ is close to the full accuracy with which we expect to be able to compute residuals in double precision and, so, going to higher orders becomes problematic if we require stricter stopping tolerances. This is similar to what we observed for the 2D case in \cite{HighOrderASMStokes}.

\begin{table}[ht]
\begin{center}
\captionsetup{justification=centering,margin=0.2cm}
\caption{{\small Convergence of FGMRES(30) preconditioned with monolithic multigrid for 3D Stokes equations with the manufactured solution given in \eqref{eq:3DManuFactSol}.}}
\label{table:PrecondSolverVankaStar3DStokes}
{
\begin{tabular}{|c|c||c|c|c||c|c|c|}
\hline

$k$   & $\ell$  &  rtol    & $\|{\bf u} - {\bf u}_h\|_1$ & $\|p-p_h\|_0$ &    Iter   &   Time    &   Memory     \\ \hline \hline
    
\multirow{2}{6pt}{3}     &  2   & $10^{ -8 } $  &   $3.90\times 10^{ -4 } $  &  $4.06\times 10^{ -4 } $ & 14   &  0.57  &   18.26 GB     \\\cline{2-8} 
  &  3 & $10^{ -8 } $    &   $5.10\times 10^{ -5 } $  &  $3.90\times 10^{ -5 } $ &  13  &  0.96   &   145.7 GB     \\\hline \hline
 
\multirow{2}{6pt}{4}      &   2  &  $10^{ -10 } $   &   $1.44\times 10^{ -5 } $  &  $1.88\times 10^{ -5 } $ &  16  &  3.70  &  58.6 GB    \\\cline{2-8} 
   &  3   &$10^{ -10 } $   &   $9.47\times 10^{ -7 } $  &  $1.03\times 10^{ -6 } $&    15  & 3.69   &  469.7 GB     \\\hline \hline

\multirow{2}{6pt}{5}       &  2   &  $10^{ -12 } $  &   $4.46\times 10^{ -7 } $  &  $5.89\times 10^{ -7 } $&  20  &  5.42  &  155.8 GB     \\\cline{2-8} 
   & 3   &   $10^{ -12 } $   &   $1.45\times 10^{ -8 } $  &  $1.56\times 10^{ -8 } $&  19  & 10.50   &   1.23 TB     \\\hline \hline
 
\multirow{2}{6pt}{6}        &  2   & $10^{ -13 } $   &   $1.23\times 10^{ -8 } $  &  $1.84\times 10^{ -8 } $&  22  &  17.66  &  358.6 GB     \\\cline{2-8} 
    &  3  &  $10^{ -13 } $   &   $2.01\times 10^{ -10 } $  &  $3.13\times 10^{ -10 } $&  20  & 20.80   &  2.9 TB     \\\hline

\end{tabular}
}

\end{center}
\end{table}

In terms of performance of the solvers, we see increasing numbers of iterations, in the column labelled ``Iter'', as we increase $k$ (and decrease the stopping tolerance), but slight decreases as we increase the number of refinements for a fixed $k$.
For fixed level of refinement, we see growing times to solution as we increase $k$, reported in the column labelled ``Time'' with data in minutes. This is due to the combination of growing numbers of DoFs with increased $k$ and growing numbers of iterations due to the decreasing stopping tolerances. Overall, we see fairly reasonable weak time scaling results from 2 to 3 levels of refinement for $k=4$ and $k=6$, but some less good weak scaling for $k=5$, which may simply be due to configuration of the machine or other factors out of our control. 
The downside to using the patch-based relaxation schemes within the solvers is the high memory costs. We report the measured memory usage for each solver in column labelled ``Memory'', measured using the slurm \texttt{seff} command.  Here, we see about a factor of 8 increase in memory requirements if we fix $k$ and increase the number of refinements, but also a strong increase in the amount of memory needed as we increase $k$ at a fixed level of refinement. 
For 3D problems, we note that the number of DoFs scales roughly like $k^3$, while the number of nonzeros in the system matrix scales like $k^6$, so expected memory increases with $k$ should be at least by factors of $4^3/3^3 \approx 2.4$ from $k=3$ to $k=4$, $5^3/4^3 \approx 2$ from $k=4$ to $k=5$ and $6^3/5^3 \approx 1.7$ from $k=5$ to $k=6$ and at most by the squares of these ratios, depending on whether they are dominated by vector or matrix storage.  The increases in measured memory usage reported in Table~\ref{table:PrecondSolverVankaStar3DStokes} are, indeed, in this range, with a measured ratio of $2.4\times$ for $\ell=3$ going from $k=5$ to $6$, which is between 1.7 and 2.9.

\subsection{Navier-Stokes Equations}

Table \ref{table:NKMVankaStar3DNavierStokes} shows similar data for the solver 
for the 3D lid-driven-cavity problem at Reynolds number of 100. In this table, we report the number of nonlinear iterations, ``NonLin'', of the outer inexact Newton solvers,
the total number of linear iterations, ``Linear'', over all nonlinear iterations, the required computational time in minutes, ``Time'', and the measured memory utilized, ``Memory''.  For the linear solves, we use the same algorithm as considered in the linear case above.  For the outer nonlinear solve, we use the same relative residual tolerances (as a function of $k$) as reported in Table~\ref{table:PrecondSolverVankaStar3DStokes} noting, of course, that we do not have an analytical solution in this case against which to check the accuracy.

We note that with $\ell = 2$, the grid likely underresolves the boundary layer at $Re = 100$, but at $\ell = 3$ and, in particular, for larger values of $k$, we expect to see fully resolved simulations.  Nonetheless, we see consistent solver performance across all values of $k$ and $\ell$, with almost no variation in the nonlinear iteration counts.
The total linear iteration counts scale very similarly to those for the (linear) Stokes case in Table \ref{table:PrecondSolverVankaStar3DStokes}. This fact shows the effectiveness of the Eisenstat-Walker method for choosing the stopping tolerance for linear solves within the inexact Newton method to avoid oversolving. 
 We also notice that the memory costs are somewhat higher in Table~\ref{table:NKMVankaStar3DNavierStokes} than those in Table~\ref{table:PrecondSolverVankaStar3DStokes}, which can be expected because we are storing more information associated with each linearization.  Nonetheless, the growth in memory usage is still generally within the bounds discussed above, albeit somewhat closer to the matrix storage bound than the vector one.  The times reported in Table~\ref{table:NKMVankaStar3DNavierStokes} are substantially larger than those reported in Table~\ref{table:PrecondSolverVankaStar3DStokes}. For instance at $(k,\ell) = (6,3)$, the time required for the Stokes case is about 20 minutes while it is more than 2 hours for the Navier-Stokes case.  This suggests that solver setup plays a more dominant role in the nonlinear case, since we have to form 6 linearization matrices and reform the preconditioner each time.  Notably, the factor-four scaling of CPU times for $\ell = 3$ going from $k=5$ to $6$ is slightly larger than the $k^6$ scaling of matrix and preconditioner assembly, suggesting significant time is spent factoring the patch and/or coarsest-grid systems in the overall solution time.  Nonetheless, achieving solution of a nontrivial nonlinear problem with over 25M DoFs in just over 2 hours on 48 cores seems acceptable if one has the available memory for this algorithm.

\begin{table}[t]
\begin{center}
\captionsetup{justification=centering,margin=0.2cm}
\caption{{\small Convergence properties of Newton-Krylov-Multigrid solver for 3D Navier-Stokes lid-driven cavity model.}}
\label{table:NKMVankaStar3DNavierStokes}
{
\begin{tabular}{|c|c||c|c|c|c|}
\hline
 $k$  &  $\ell$     &    NonLin   &   Linear  &   Time    &   Memory  \\ \hline \hline
    
\multirow{2}{6pt}{3}    &   2   &  6& 15     & 1.1  &  20.66 GB   \\\cline{2-6} 
  &   3   &  5&  14    &  1.8   & 168.39 GB \\\hline \hline
 
\multirow{2}{6pt}{4}     &   2   & 6&  18    &  4.9   & 78.82 GB   \\\cline{2-6} 
   &    3   & 6&  19     & 8.8   & 661.01  GB  \\\hline \hline

\multirow{2}{6pt}{5}      &  2   &  6& 21     &  17.6    &  242.14 GB   \\\cline{2-6} 
     &  3   &  6&  21    &  34.7    & 2.03 TB  \\\hline \hline
 
\multirow{2}{6pt}{6}    &   2   &  6&  24    &  100.4    & 633.59  GB  \\\cline{2-6} 
   &   3   &  6&  24    &  135.3    & 5.23 TB \\\hline

\end{tabular}
}

\end{center}
\end{table}


\section{Conclusion and future work}\label{sec:conclusion}

In this work, we study the extension of monolithic multigrid methods with subspace decomposition (Vanka) relaxation for higher-order Taylor-Hood finite-element discretizations of the Stokes and Navier-Stokes equations.  The key ingredient introduced is a new patch construction for the Vanka relaxation, named Vanka-star relaxation, that introduces overlap into the pressure variables in the subspace decomposition.  We show numerically that this leads to both improved iteration counts and time-to-solution for a range of discretization orders, and that the improved linear solver performance leads to greatly improved nonlinear solver performance for the Navier-Stokes equations.  This improvement enables simulation of three-dimensional flows demonstrating the efficiency of the approach.

Natural opportunities for future work include extending this relaxation scheme to more complicated flow scenarios, including more challenging three-dimensional flows, as well as flows of complex fluids, such as thermal flows or magnetohydrodynamics.  Two key algorithmic innovations that are needed are improvements in forming the patch systems, where it is possible that surrogate models can be used to avoid forming all of the patch systems, as was investigated in~\cite{harper2023compressionreducedrepresentationtechniques}.  A further question in this direction is whether inexact or iterative solves can be used in place of the dense LU factorizations considered here, similar to the ``diagonal Vanka'' techniques considered in~\cite{NavierStokesVolkerTobiska} and elsewhere.

\section*{Acknowledgments}
The work of SM was partially supported by an NSERC Discovery Grant.  This research was enabled in part by ACENET (\url{ace-net.ca}) and the Digital Research Alliance of Canada (\url{alliancecan.ca}).

\end{document}